\documentclass[a4paper,10pt]{article}

\usepackage{float}
\usepackage{amsmath,amssymb,amsthm,amscd,algorithmic,algorithm}
\usepackage[mathscr]{eucal}
\usepackage[all]{xy}
\usepackage{graphicx}
\usepackage{appendix}
\usepackage{array}
\usepackage{booktabs}
\usepackage{color}
\usepackage{hyperref}
\hypersetup{
colorlinks=true,
linkcolor=blue,
citecolor=blue,
urlcolor=blue
}

\makeatletter
    
    \@addtoreset{equation}{section}
  \makeatother

\newcommand{\mcal}{\mathcal}
\newcommand{\mbf}{\mathbf}

\newcommand{\mbb}{\mathbb}
\newcommand{\mrm}{\mathrm}
\newcommand{\vphi}{\varphi}

\newcommand{\aet}{\mathrm{\acute{e}t}}
\newcommand{\Q}{\mathbb{Q}}
\newcommand{\Z}{\mathbb{Z}}
\newcommand{\tor}{\mathrm{tor}}
\newcommand{\cO}{\mathcal{O}}

\newcommand{\OO}{\mathcal{O}}
\newcommand{\pp}{\mathfrak{p}}
\newcommand{\bb}{\mathfrak{b}}
\newcommand{\mm}{\mathfrak{m}}
\newcommand{\minfty}{\mathfrak{m}_{\infty}}
\newcommand{\Ah}{\widehat{A}}
\newcommand{\Ab}{\overline{A}}

\newcommand{\Fr}{\mathrm{Fr}}
\DeclareMathOperator{\Gal}{Gal}
\DeclareMathOperator{\Art}{Art}
\DeclareMathOperator{\Nr}{Nr}

\DeclareMathOperator{\GL}{GL}
\DeclareMathOperator{\red}{red}

\DeclareMathOperator{\rank}{rank}

\DeclareMathOperator{\Cl}{Cl}

\newcommand{\Up}[2]{U_{#1}^{(#2)}}

\newtheorem{theorem}{Theorem}[section]
\newtheorem{corollary}[theorem]{Corollary}
\newtheorem{lemma}[theorem]{Lemma}
\newtheorem{proposition}[theorem]{Proposition}

\theoremstyle{definition}
\newtheorem{definition}[theorem]{Definition}
\newtheorem{remark}[theorem]{Remark}

\newtheorem{example}[theorem]{Example}

\newtheorem*{acknowledgments}{Acknowledgments}

\counterwithin{table}{section}
\counterwithin{algorithm}{section}

\title{
Local decomposition of torsion on abelian varieties and its global applications}

\author{Yoshiyasu Ozeki\footnote{
Faculty of Science, Kanagawa University,
3-27-1 Rokkakubashi, Kanagawa-ku, Yokohama-shi, Kanagawa 221-8686, JAPAN
\endgraf
e-mail: {\tt ozeki@kanagawa-u.ac.jp}
}  
and Manabu Yoshida\footnote{
Faculty of Mathematics Education, Osaka Kyoiku University, 
4-698-1 Asahigaoka, Kashiwara, Osaka 582-8582, JAPAN
email: {\tt yoshida-m95@cc.osaka-kyoiku.ac.jp}
\endgraf
}}

\begin{document}
\maketitle

\begin{abstract}
In this paper, we establish a decomposition theorem for torsion of abelian varieties over Galois extensions of complete discretely valued fields of mixed characteristic. Under mild potential semistability and pro-nilpotent inertia assumptions, the torsion splits into a ramified $p$-primary part, an unramified $p$-primary part, and a tame prime-to-$p$ part. We call this result the local decomposition theorem.
We describe the ramified part using Raynaud uniformization and, under further assumptions, determine its group structure and obtain bounds for its order.

We give two main global applications of our local decomposition theorem. First, for any $g$-dimensional abelian variety $A$ over $\mathbb{Q}$, we prove that
$A(\mathbb{Q}_{\infty,p})_{\mathrm{tor}}=A(\mathbb{Q})_{\mathrm{tor}}$
whenever $p$ exceeds an explicit constant depending only on $g$, where $\mathbb{Q}_{\infty,p}$ denotes the cyclotomic $\mathbb{Z}_p$-extension of $\mathbb{Q}$. Second, we completely determine the group structure of $J_0(N)(\mathbb{Q}^{\mathrm{ab}})_{\mathrm{tor}}$ for several prime levels $N$.
We also obtain explicit bounds for torsion growth and finite-layer stabilization results over relative ray class field towers.

\end{abstract}

    \tableofcontents


\section{Introduction}

Let $F$ be either a number field or a $p$-adic field, and let $L/F$
be an algebraic extension, typically of infinite degree.
Let $A$ be an abelian variety over $F$.
A fundamental problem is to study the finiteness of
$A(L)_{\mathrm{tor}}$ and, when it is finite, to determine its group
structure and obtain bounds for its order.

There are several precedents for studying torsion over infinite
extensions. For elliptic curves over $\mathbb Q$, Chou \cite{Cho19}
classifies the groups
$E(\mathbb Q(\mu_\infty))_{\mathrm{tor}}$.
For elliptic curves over $\mbb{Q}_p$ with good reduction,
the authors \cite{OzYo26} determine the groups
$E(\mbb{Q}_p)_{\mathrm{tor}}$ and
$E(\mbb{Q}_p(\mu_{p^\infty}))_{\mathrm{tor}}$.
Here $\mu_\infty$ denotes all roots of unity, and
$\mu_{p^\infty}$ denotes the $p$-power roots of unity.
In arbitrary dimension, Ribet's theorem \cite[Appendix]{KaLa81}
shows that
$A(F(\mu_\infty))_{\mathrm{tor}}$
is finite for every abelian variety $A$ over a number field $F$,
while Imai \cite{Ima75} proves that
$A(K(\mu_{p^\infty}))_{\mathrm{tor}}$
is finite for an abelian variety $A/K$ with potentially good reduction,
where $K$ is a $p$-adic field.
Explicit upper bounds for
$\#A(K)_{\mathrm{tor}}$
are given by Clark--Xarles \cite{ClXa08}
when $A/K$ has anisotropic reduction,
a condition including potentially good reduction.
The first author \cite{Oze24-2} obtains explicit bounds for
torsion of CM abelian varieties over Lubin--Tate extensions.

\medskip 
We now introduce the main structural result of this paper.
Our approach is to separate the torsion acquired over an extension
into a ramified $p$-primary part, an unramified $p$-primary part,
and a tame prime-to-$p$ part.
For this purpose, we introduce a new method, which we call
\emph{local decomposition}.
A notable feature of this approach is that it can be applied to abelian
varieties of {\it arbitrary dimension}.

Throughout this paper,  let $p>2$ be an odd prime.
Let $K$ be a complete discretely valued field of mixed characteristic
$(0,p)$ with perfect residue field.
For every finite extension $K'/K$, we denote by $e_{K'}$
the absolute ramification index of $K'$.
We also denote by $K^{\mrm{ur}}$ and $K^{\mrm{tame}}$
the maximal unramified extension and the maximal tamely ramified extension
of $K$, respectively.
For any abelian variety $A_{/K}$ and any Galois extension $L/K$, we set 
$$
D_{A,K}(L):=\sum_{\sigma\in I(L/K)} (\sigma -1)A(L)[p^{\infty}],
$$
which is a $G_K$-stable submodule of $A(L)[p^{\infty}]$. 
We will show in Section \ref{LD_st} that,
under certain assumptions, $D_{A,K}(L)$ admits a natural interpretation
in terms of rigid analytic uniformization
(see Proposition \ref{LD:semistable} (1) and Corollary \ref{rigid_D}).
For example, suppose that $e_K<p-1$ and that $L/K$ is a Galois extension
with pro-nilpotent inertia.
This includes, in particular, the cases where $L/K$ is abelian, 
tamely ramified or pro-$p$.
Then, $D_{A,K}(L)$ coincides with the
$p$-primary torsion in the reduction kernel
when $A$ has good reduction.
More generally, if $A$ has semistable reduction,
then $D_{A,K}(L)$ is identified, via rigid analytic uniformization, 
with the $p$-primary torsion
in the reduction kernel of the associated
semi-abelian variety.
$D_{A,K}(L)$ is isomorphic to $\mu_{p^\infty}(L)^g$
when $A$ has split totally degenerate reduction
of dimension $g$;
see Corollary \ref{rigid_D}.

We clearly have a natural inclusion
$$
A(L)_{\mrm{tor}}\supset 
D_{A,K}(L)
+A(L\cap K^{\mrm{ur}})[p^{\infty}]
+A(L\cap K^{\mrm{tame}})_{p'},
$$
where the subscript $p'$ denotes the prime-to-$p$ part.
When the inclusion above is an equality and the sum on the right-hand side is direct, 
we call the resulting decomposition
the local decomposition of $A(L)_{\mathrm{tor}}$.

\begin{theorem}\label{thm:intro-local}
Let $K$ be a complete discretely valued field of mixed characteristic
$(0,p)$ with perfect residue field.
Let $L/K$ be a Galois extension with pro-nilpotent inertia.
Let $A$ be an abelian variety over $K$ 
which acquires  semistable reduction
over a finite extension $K'/K$ satisfying $e_{K'}< p-1$.
Then, 
\[
A(L)_{\mrm{tor}}
= 
D_{A,K}(L)
\oplus A(L\cap K^{\mrm{ur}})[p^{\infty}]
\oplus A(L\cap K^{\mrm{tame}})_{p'}.
\]
Moreover, for every $I(L/K)$-stable submodule
$N\subset A(L)[p^\infty]$,
\[
 N=\bigl(N\cap D_{A,K}(L)\bigr)\oplus N^{I(L/K)}.
\]
\end{theorem}


The situation becomes particularly explicit when  
$L/K$ is a totally ramified abelian extension.
In this case, under some mild assumptions, 
the group structure of $D_{A,K}(L)$ can be determined completely, and
an explicit upper bound for its order can also be obtained.
See Theorem \ref{FG:Main} for the precise statement.
The most typical example of Theorem \ref{thm:intro-local} 
is the cyclotomic extension.
Let $K/{\mbb{Q}_p}$ be a finite unramified extension with residue field $k$, and 
let $A/K$ have good reduction.
Local decomposition gives a non-canonical isomorphism of finite abelian groups
\[
 A(K(\mu_{p^\infty}))_{\mathrm{tor}}
 \simeq A(K)_{\mathrm{tor}}
       \oplus\overline A^{\vee}(k)[p^\infty]
\]
(see Theorem \ref{hatG_LD}).

\medskip 

Let us give some applications of local decompositions.
For elliptic curves $E$ over $\mathbb Q$, Chou--Daniels--Krijan--Najman 
proved in \cite[Theorem 1.1]{ChDaKrNa21} that
$
E(\mathbb{Q}_{\infty,p})_{\rm tor}=E(\mathbb Q)_{\rm tor}
$
for every $p\ge5$. Here, $\mbb{Q}_{\infty,p}/\mbb{Q}$ is the cyclotomic $\mbb{Z}_p$-extension. 
Our local decomposition yields the following higher-dimensional analogue
for abelian varieties over number fields.
For any $g\ge 1$, put 
$$
m_g=\mrm{max}\left\{\mrm{lcm}(n_1,\dots ,n_r)
\mid r>0,\ n_1,\dots, n_r>0,\  \sum^r_{i=1}\varphi(n_i)\le 2g\right\}.
$$
For example, $m_1=6$, $m_2=12$, $m_3=30$ and $m_4=60$.
See Table~\ref{tab:mg} for further values of $m_g$.
\begin{theorem}[Corollary \ref{cor:cyclotomic-rigidity}]
\label{cor:cyclotomic-rigidity_intro}
Let $F$ be a number field and put
$
F_{\infty,p}:=F\mathbb Q_{\infty,p}.
$
Assume that there exists a finite  place $v$ of $F$ above $p$
which is unramified over $\mathbb Q$.
Let $A_{/F}$ be a $g$-dimensional abelian variety.
\begin{itemize}
\item[{\rm (1)}]
If $A$ has semistable reduction at $v$, then
$$
A(F_{\infty,p})_{\rm tor}=A(F)_{\rm tor}.
$$

\item[{\rm (2)}]
If $p>m_g+1$, then
$$
A(F_{\infty,p})_{\rm tor}=A(F)_{\rm tor}.
$$
\end{itemize}
\end{theorem}

As illustrated by Theorem 1.2,
the local decomposition theorem is particularly useful
for reducing the study of torsion over infinite extensions
to that over smaller fields.
Combined with the results of Mazur~\cite{Maz77},
it also enables us to completely determine the torsion 
subgroups of certain modular Jacobians over infinite extensions.
Let $N$ be a prime and put $J=J_0(N)$.
We determine the group structure of $J(E)_{\mrm{tor}}$
in the following cases:
\begin{itemize}
\item[{\rm (a)}]
(Example~\ref{Jac_Ex})
$E$ is the maximal pro-$7$ Galois extension of
$\mbb Q(\sqrt{-1})$ unramified outside $7$, and
$N\in\{23,29,31,37\}$.
\item[{\rm (b)}]
(Table~\ref{tab:jacobian-torsion})
$E=\mbb Q^{\mrm{ab}}$, and
$N\in\{11,19,23,29,37,47,53,59,61,71,79,83\}$.
\end{itemize}
Precisely, writing $n$ for the numerator of $(N-1)/12$,
we have
$$
J(E)_{\mrm{tor}}\simeq
\begin{cases}
\mbb Z/n\mbb Z & \text{in case (a)},\\
(\mbb Z/n\mbb Z)^{2} & \text{in case (b)}.
\end{cases}
$$
In case (a), by a result of Mazur and Rubin~\cite{MaRu08},
the Galois group of $E/\mbb Q(\sqrt{-1})$
is a free pro-$7$ group on two generators
and is thus non-abelian.
Here we focus on the two-dimensional Jacobians
as particularly accessible examples.
In all these cases, $J(E)_{\mrm{tor}}$
coincides with the cuspidal subgroup of $J$.
Further examples can be obtained
by the same method together with computations in MAGMA.

In case (b), the theory of the Eisenstein ideal,
developed by Mazur and Ribet,
plays a central role in the study of torsion
over $\mbb Q^{\mrm{ab}}$.
Our approach is different in that it combines
the results of Mazur~\cite{Maz77}
with the local decomposition theorem
without relying on the structure theory
of the Eisenstein ideal.
To the best of our knowledge, no explicit determination
of the entire torsion subgroup over $\mbb Q^{\mrm{ab}}$
for the higher-dimensional cases listed in (b)
has appeared in a published paper.
After obtaining these results, we learned of
a talk announcement by Ribet,
in which he states that, for any prime level $N$,
the torsion subgroup of $J_0(N)$ over
$\mbb Q^{\mrm{ab}}$ coincides with the kernel
of the Eisenstein ideal.
Ribet kindly confirmed this result
in private correspondence.
It remains an interesting problem to investigate
to what extent the local decomposition theorem
can be applied to the case of composite levels.
\medskip

Finally, we apply local decomposition and the description of formal
torsion to ray class field towers.
The following special case gives an explicit bound for torsion growth
and a stabilization level, both expressed in terms of a single value
of the local Frobenius polynomial.
\begin{theorem}\label{thm:intro-ray}
Let $F$ be a number field, $p>2$ a prime, and $\pp=(\pi)$ a prime of $F$ above $p$ 
with $\pi$ totally positive. 
Assume that $K=F_\pp$ is unramified over $\Q_p$ and that $A/F$ is a $g$-dimensional 
abelian variety with good ordinary reduction at $\pp$. 
For $n\ge0$, let $F_n$ be the ray class field of modulus $\pp^n \mm_\infty$, 
where $\mm_\infty$ is the product of the real places of $F$, 
and put $H^+=F_0$, the narrow Hilbert class field, and $F_\infty=\bigcup_{n\ge0}F_n$.

Write $P_{A,K}$ for the Frobenius polynomial of the reduction $\overline A$ 
over the residue field $k_K$ of $K$, and set $s=\mrm{Nr}_{K/\Q_p}(\pi)$. 
If $P_{A,K}(s) \not=0$, then $m:=v_p(P_{A,K}(s))-g [K:\Q_p] \in \Z_{\ge0}$ 
satisfies
\[
[A(F_\infty)_\tor:A(H^+)_\mrm{tor}]\le p^m,\qquad A(F_\infty)_\tor=A(F_m)_\tor,
\]
where $v_p(p)=1$. 
If $s=\# k_K$, the nonvanishing condition is automatic and 
$m=v_p(\# \overline{A}(k_K))$. 
\end{theorem}

More general results for relative ray class
towers, including the semistable case and
descriptions via Smith invariants,
are given in Section \ref{sec:ray}.

\medskip

The paper is organized as follows.
Section~2 establishes the local decomposition theorem,
including the semistable case.
Section~3 studies formal torsion over totally ramified
extensions and gives explicit descriptions and bounds.
Section~4 presents global applications, including
rigidity results over pro-$p$ extensions and
torsion of modular Jacobians.
Finally, Section~5 treats torsion over ray class field towers.

\begin{acknowledgments}
The authors are grateful to Takashi Hara for his valuable advice on connections 
with Iwasawa theory at an early stage of this work. 
They also thank Kenneth A. Ribet for kindly sharing his lecture slides 
containing an outline of the proof of his result on cyclotomic torsion points of modular Jacobians.
\end{acknowledgments}

\vspace{5mm}

\noindent
{\bf AI disclosure.}\ 
In this paper, we used ChatGPT for English proofreading, 
organizing notation and mathematical arguments, and checking proofs. 
The central idea of the local decomposition theorem arose from the authors' 
investigations in their previous paper \cite{OzYo26} and was developed 
and proved by the authors themselves.
For Section \ref{Sec_FG_ur}, we discussed technical details concerning 
finite flat group schemes with ChatGPT. 
Corollary \ref{cor:cyclotomic-rigidity} was originally established by the authors, 
but discussions with ChatGPT helped us arrive at its present 
formulation in terms of $m_g$. 
For Section \ref{Sec_Jac}, ChatGPT provided an initial 
framework for the case $N=23$, which the authors subsequently 
refined and generalized to a wider range of prime levels $N$.
The authors take full responsibility for the mathematical content of this paper.

\section{Local decomposition theorems}
\label{secton:LD}
Let $p$ be a prime.
Let $K$ be a complete discretely valued field of mixed characteristic $(0,p)$
with perfect residue field $k$.
For any finite extension $K'/K$, 
we denote by $e_{K'}$ the absolute ramification index of $K'$.
Let $g>0$ be an integer.
Let $A$ be a $g$-dimensional  abelian variety over $K$. 
The goal of this section is to explain the analysis of $D_{A,K}(L)$
for semistable reduction cases
(cf.\ Proposition \ref{LD:semistable} (1) and Corollary \ref{rigid_D})
and to prove Theorem \ref{thm:intro-local}.

\subsection{The case of semistable reduction}
\label{LD_st}
In this section, we assume that {\it $A_{/K}$ has semistable reduction}.
Consider the non-Archimedean rigid uniformization of $A$; 
there exist the following data (cf.\  \cite[Definition 1.1 and Theorem 1.2]{BoXa96}).  
\begin{itemize} 
\item[(i)] $G$ is a semi-abelian variety of dimension $g$ that 
fits into an exact sequence
of $K$-group schemes
$0\to T \to G \to B \to 0$ where $T$ is a torus of rank $m$ and $B$
is an abelian variety which has potentially good reduction,
\item[(ii)] a closed immersion of rigid $K$-groups $N^{\rm an}\hookrightarrow G^{\rm an}$ 
where  $N$ is a group scheme which is isomorphic to $\mbb{Z}^{\oplus m}$
 after a finite base extension. Here, the superscript  
 ``an'' denotes rigid analytification, and  
\item[(iii)] a faithfully flat morphism $u\colon G^{\rm an}\to A^{\rm an}$ of rigid $K$-groups 
with kernel $N^{\rm an}$.
\end{itemize}
Since $A$ has semistable reduction over $K$,
it follows from \cite[Proposition 4.1]{BoLi99} that 
$B$ has good reduction over $K$, $T$ splits over a finite unramified extension $K'/K$
and $N$ is constant over $K'$.
For any Galois extension $L/K$ with residue field $k_L$, we set 
$$
\widehat{G}(L):=\mrm{ker}(G(L)\to B(L)\to \bar{B}(k_L)),
$$
the kernel of the reduction map of $G(L)$. This is a $G_K$-stable submodule of $G(\overline{K})$
and, by definition, we have 
$\widehat{G}(\overline{K})^{G_L}=\widehat{G}(L)$.
We have exact sequences
\begin{equation}
\label{rigid_ex}
0\to N(\overline{K}) \to G(\overline{K})\to A(\overline{K})\to 0 \quad \mrm{and} \quad 
0\to T(\overline{K}) \to G(\overline{K})\to B(\overline{K})\to 0
\end{equation}
of $G_K$-modules (since the morphisms
$N^{\rm an}\to G^{\rm an}\to A^{\rm an}$
are rigid analytic, the exactness of the former sequence
might not be well known; see the proof of
\cite[Theorem 2.3]{ClXa08} for this).

\begin{example}
\label{ex_D}
\begin{itemize}
\item[(1)] If $A_{/K}$ has good reduction,
then $G=B=A$, $T=0$ and $N=0$.
Thus $\widehat{G}(L)$ 
coincides with the kernel of the reduction map on $A(L)$.
\item[(2)] If $A$ is a $g$-dimensional abelian variety with split totally
degenerate reduction, then $G=T=\mbb{G}_m^g$, $B=0$ and $N=\mbb{Z}^g$.
Thus $\widehat{G}(L)\simeq (L^{\times})^g$ as $G_K$-modules.
\end{itemize}
\end{example}

\medskip

We consider the following condition $(C1)$ for $A_{/K}$ and $L/K$.

\begin{itemize}
\item[{\bf (C1):}] 
There exists a cofinal family $\mathcal{C}$ of finite Galois subextensions
$L'/K$ of $L/K$ such that, for every $L'\in\mathcal{C}$, 
one can choose $z\in Z(I(L'/K))$
so that 
$z-1$ acts bijectively on $\widehat{G}(L')[p]$.
\end{itemize}
Here, cofinal means that every finite subextension of $L/K$
is contained in some $L'\in\mathcal{C}$. 

\begin{proposition}\label{LD:semistable}
Let the notation be as above. Assume that {\rm (C1)} is satisfied for $A_{/K}$ and $L/K$.

\begin{itemize}
\item[{\rm (1)}] We have $D_{A,K}(L)=u(\widehat{G}(L)[p^{\infty}])$.
\item[{\rm (2)}] For every $I(L/K)$-stable submodule
$N\subset A(L)[p^\infty]$, 
\[
 N=\bigl(N\cap u(\widehat{G}(L)[p^{\infty}])\bigr)\oplus N^{I(L/K)}.
\]
\end{itemize}
\end{proposition}

Before the proof, we give some preliminaries.

\begin{itemize}
\item[--] It follows from \eqref{rigid_ex}
that, for any $0<n\le \infty$, we have an exact sequence 
\begin{equation}
\label{rigid_ex1}
0\to G(\overline{K})[p^n]\overset{u}{\rightarrow} A(\overline{K})[p^n]\to 
N(\overline{K})\otimes_{\mbb{Z}}p^{-n}\mbb{Z}/\mbb{Z}      
\end{equation}
of $G_K$-modules and $N(\overline{K})$ is unramified.
(Here, we set $p^{-\infty}\mbb{Z}/\mbb{Z}:=\mbb{Z}[1/p]/\mbb{Z}$.)
\item[--] By definition, 
for any $0<n\le \infty$, 
we have an exact sequence 
\begin{equation}\label{rigid_ex2}
0\to T(\overline{K})[p^n]\to \widehat{G}(\overline{K})[p^n]
\overset{v}{\rightarrow} \widehat{B}(\overline{K})[p^n]     
\end{equation}
of $G_K$-modules and $T(\overline{K})[p^n]\simeq \mu_{p^n}(\overline{K})^m$ as $I_K$-modules.
Here, $\widehat{B}(\overline{K})$ is the kernel of the reduction map on 
$B(\overline{K})$.
\end{itemize}

\begin{lemma}\label{hatG1}
Let $A_{/K}$ be an abelian variety with semistable reduction, 
and retain the notation  above.
Let $L/K$ be a Galois extension and 
$0<s\le \infty$.
\begin{itemize}
\item[{\rm (1)}] We have 
$u(\widehat{G}(\overline{K})[p^s])^{G_L}
=u(\widehat{G}(L)[p^s])$.
\item[{\rm (2)}] We have 
$u(\widehat{G}(L)[p^s])[p^t]
=u(\widehat{G}(L)[p^t])$
for any integer $0<t\le s$.
\end{itemize}
\end{lemma}
\begin{proof}
The result follows from the injectivity of the  $G_K$-equivariant homomorphism 
$u\colon G(\overline{K})[p^s]\to A(\overline{K})[p^s]$.
\end{proof}

\begin{lemma}\label{hatG2}
Let $A_{/K}$ be an abelian variety with semistable reduction, 
and retain the notation above.
Let $n>0$ be an integer and 
$L/K$ a finite Galois extension.
Then 
$$
(\sigma -1)A(L)[p^n]\subset u(\hat{G}(L)[p^n])
$$
for every $\sigma \in I(L/K)$.
\end{lemma}

\begin{proof}
Let $P\in A(L)[p^n]$ and 
take any $\sigma\in I(L/K)$.
Choose a lift $\hat{\sigma}\in I_K$ of $\sigma$.
Since $I_K$ acts trivially on $N(\overline{K})$,
the element $\sigma(P)-P=\hat{\sigma}(P)-P$ belongs to the kernel of 
$A(L)[p^n]\to N(\overline{K})/p^nN(\overline{K})$.
This implies 
\begin{equation}
\label{lem_ex1}
\sigma(P)-P\in u(G(L)[p^n])
\end{equation}
by \eqref{rigid_ex1}.

Choose a lift $\hat{P}$ of $P$ under the surjection
$u\colon G(\overline{K})\to A(\overline{K})$.
We have $p^n\hat{P}\in \ker u =N(\overline{K})$. 
Since $I_K$ acts trivially on $N(\overline{K})$,
$p^n(\hat{\sigma}\hat{P}-\hat{P})=0$,
and hence 
$\hat{\sigma}\hat{P}-\hat{P}\in G(\overline{K})[p^n]$.
On the other hand, since $I_K$ acts trivially 
on $\bar{B}(\bar k)$, 
we see $\hat{\sigma}\hat{P}-\hat{P}\in \widehat{G}(\overline{K})$.
Therefore, 
$
\hat{\sigma}\hat{P}-\hat{P}\in \widehat{G}(\overline{K})[p^n],
$
and thus 
\begin{equation}\label{lem_ex2}
\sigma(P)-P=u(\hat{\sigma}(\hat{P})-\hat{P})
\in u(\widehat{G}(\overline{K})[p^n]).
\end{equation}
By \eqref{lem_ex1} and \eqref{lem_ex2},
we obtain 
\begin{align*}
(\sigma-1)(P)
&\in u(G(L)[p^n])\cap u(\widehat{G}(\overline{K})[p^n])\\
&\subset u(\widehat{G}(\overline{K})[p^n])^{G_L}
=u(\widehat{G}(L)[p^n]),
\end{align*}
where the last equality follows from Lemma \ref{hatG1} (1). 
\end{proof}

\begin{lemma}   
\label{lemmaA}
Let $D$ be a finite module of $p$-power order. 
Let $u$ be an endomorphism of $D$.
Then, $u$ is bijective on $D$ if and only if 
its restriction to $D[p]$ is bijective.
\end{lemma}

\begin{proof}
Suppose that $u$ is bijective on $D[p]$. 
Assume that $u$ is not injective on $D$.
Then there exists a non-zero $x\in \ker (u\ \mrm{on}\ D)$.
If we denote by $p^m$ the order of $x$, then 
$0\not=p^{m-1}x\in D[p]\cap \ker (u\ \mrm{on}\ D)$,
which contradicts the injectivity of $u$ on $D[p]$. 
Thus $u$ is injective on $D$, and hence  bijective since $D$ is finite.
\end{proof}

\begin{lemma}
\label{lemmaB}
Let $\Gamma$ be a finite group and $M$ a finite $\Gamma$-module of $p$-power order.
Let $D$ be a $\Gamma$-stable submodule of $M$. 
Assume the following.
\begin{itemize}
\item[{\rm (i)}] $(\sigma-1)M\subset D$ for every $\sigma \in \Gamma$, and 
\item[{\rm (ii)}] there exists an element $z\in Z(\Gamma)$ such that 
$z-1$ acts bijectively on $D[p]$.
\end{itemize}
Fix an element $z$ as in {\rm (ii)}. 
Then the following hold.
\begin{itemize}
\item[{\rm (1)}] $w:=z-1\colon D\to D$ is bijective.
\item[{\rm (2)}] Let $e:=w^{-1}\circ (z-1)\colon M\to D$, which is well-defined by {\rm (1)} and {\rm (i)}. 
Then, $e$ is $\Gamma$-equivariant 
and satisfies $e|_D=\mrm{id}_D$, $e^2=e$ and $\ker e=M^{\Gamma}$. 
In particular, the natural inclusion $M\supset D+M^{\Gamma}$
is an equality and the sum is direct;
\begin{equation}
\label{lemmaB_dec}
M= D\oplus M^{\Gamma}.
\end{equation}
\item[{\rm (3)}] If $N$ is a $\Gamma$-stable submodule of $M$, then
we have $e(N)=N\cap D$ and 
the restriction of the decomposition in \eqref{lemmaB_dec} to $N$ yields 
$$
N=(N\cap D)\oplus N^{\Gamma}.
$$
\end{itemize}
\end{lemma}
\begin{proof}
The assertion (1) follows from Lemma \ref{lemmaA} and (ii).
We show (2).
Since $z$ is an element of the center of $\Gamma$, the map $w=z-1$
commutes with the action of $\Gamma$, and hence so does $w^{-1}$.
Therefore, $e$ is $\Gamma$-equivariant. By definition we have $e|_D=\mrm{id}_D$.
For any $x\in M$, since $y:=e(x)\in D$,
we have 
$e^2(x)=e(y)=y=e(x)$. Thus $e^2=e$.
We next show $\ker e=M^{\Gamma}$.
The inclusion $M^{\Gamma}\subset \ker e$ is immediate from the definition of $e$.
Take any $x\in \ker e$. 
For any $\sigma \in \Gamma$, by  $(\sigma-1)(x)\in D$, $e|_D=\mrm{id}_D$ and $e$ is 
$\Gamma$-equivariant,
we have 
$$
(\sigma-1)(x)=e((\sigma-1)(x))
=(\sigma-1)(e(x))=0.
$$ 
Hence $x\in M^{\Gamma}$, and therefore 
$\ker e=M^{\Gamma}$.
For any $x\in M$, we have
$x=e(x)+(x-e(x))$, 
where $e(x)\in D$ and
$
e(x-e(x))=e(x)-e^2(x)=0.
$
Thus $x-e(x)\in\ker e=M^{\Gamma}$.
Moreover,
$D\cap M^{\Gamma}=D\cap\ker e=0$
by $e|_D=\mrm{id}_D$.
Therefore, we have 
$M=D\oplus M^{\Gamma}$
as desired.

Finally, we show (3).
Since $D$ and $N$ are $\Gamma$-stable,
we have $w(N\cap D)\subset N\cap D$.
Thus $w$ acts bijectively on $N\cap D$
since $N\cap D$ is finite.
This shows $w^{-1}(N\cap D)=N\cap D$.
By (i) and the $\Gamma$-stability of $N$, we have $(z-1)(N)\subset N\cap D$.
Hence 
$$
e(N)=w^{-1}(z-1)(N)\subset N\cap D.
$$
By $e|_D=\mrm{id}_D$, we obtain $e(N)=N\cap D$.
For any $x\in N$, we have 
$$
x-e(x)\in N\cap \ker e
=N\cap M^{\Gamma}=N^{\Gamma}.
$$
Thus the remaining assertion follows by the same argument as in the last part of the proof of (2).
\end{proof}

\begin{proof}[{Proof of Proposition \ref{LD:semistable} }]
Let $\mcal{C}$ be the cofinal family appearing in condition (C1) for $A_{/K}$ and $L/K$.
First we show (2).
We begin with the following claim.

\medskip

\noindent
{\bf Claim.}  For  any integer $n>0$, 
any $L'\in \mcal{C}$
and any $I(L'/K)$-stable submodule $N'\subset A(L')[p^n]$, 
$$
N'=(N'\cap u(\hat{G}(L')[p^n]))\oplus N'^{I(L'/K)}.
$$

\medskip

\noindent
Indeed, put  $M=A(L')[p^n]$, 
$D=u(\hat{G}(L')[p^n])$ 
and $\Gamma=I(L'/K)$.
It suffices to verify the assumptions of Lemma \ref{lemmaB}.
Assumption (i) follows from Lemma \ref{hatG2}, 
while assumption (ii) follows from condition (C1).
Thus the claim follows from Lemma \ref{lemmaB} (3).

Now let $N$ be as in (2).
For any $n>0$, put 
$$
N_n=N[p^n]=N\cap A(L)[p^n]. 
$$
Since $N_n$ is finite and $\mcal{C}$ is cofinal, we may choose 
$L_n\in \mcal{C}$ such that 
$N_n\subset A(L_n)$.  
The action of $I(L/K)$ on $N_n$ factors through $I(L_n/K)$,
and $N_n$ is an $I(L_n/K)$-stable submodule of $A(L_n)[p^n]$.
Hence  the claim gives  
\begin{equation}
\label{Nn1}
N_n=(N_n\cap u(\hat{G}(L_n)[p^n]))\oplus N_n^{I(L_n/K)}.
\end{equation}
Take any $P\in N$ and suppose $p^nP=0$ for some $n>0$.
By \eqref{Nn1}, we have 
$$
P\in N_n\subset (N\cap u(\hat{G}(L)[p^{\infty}]))+ N^{I(L/K)}.
$$
Hence 
$$
N=(N\cap u(\hat{G}(L)[p^{\infty}]))+ N^{I(L/K)}.
$$
It remains to show that the sum is direct.
Take
$$
P\in (N\cap u(\hat{G}(L)[p^{\infty}]))\cap N^{I(L/K)}
=u(\hat{G}(L)[p^{\infty}])\cap N^{I(L/K)}
$$ 
and choose $n>0$ such that $p^nP=0$. 
By Lemma \ref{hatG1} (2), $P\in u(\hat{G}(L)[p^n])$.
Since $P$ is $L_n$-rational,
it follows from Lemma \ref{hatG1} (1)
that 
$P\in u(\hat{G}(L)[p^n])^{G_{L_n}}
=u(\hat{G}(L_n)[p^n])$.
Hence $P\in u(\hat{G}(L_n)[p^n])\cap N_n^{I(L_n/K)}$.
By \eqref{Nn1}, we obtain $P=0$. This proves (2).

Next we show (1).
Let $n>0$ be an integer and  $L'\in \mcal{C}$. 
Put 
$$
D_{A,K}(L')_n=\sum_{\sigma\in I(L'/K)} (\sigma-1)A(L')[p^n]
=\sum_{\sigma\in I_K} (\sigma-1)A(L')[p^n]
$$
where the last sum is understood as a finite sum.
Let us prove 
\begin{equation}\label{D_AK}
D_{A,K}(L')_n=u(\hat{G}(L')[p^n]).
\end{equation}
The inclusion 
$D_{A,K}(L')_n\subset u(\hat{G}(L')[p^n])$
follows from Lemma \ref{hatG2}.
By condition (C1), there exists $z\in Z(I(L'/K))$
such that $z-1$ acts bijectively on $\hat{G}(L')[p]$.
Then $z-1$ acts bijectively on $u(\hat{G}(L')[p^n])$
by Lemma \ref{lemmaA}.
Therefore, 
$$
u(\hat{G}(L')[p^n])=(z-1)u(\hat{G}(L')[p^n])
\subset (z-1)A(L')[p^n]\subset D_{A,K}(L')_n.
$$
This proves \eqref{D_AK}.
By taking the union of \eqref{D_AK} over $n>0$, 
we obtain 
$$
D_{A,K}(L')=u(\hat{G}(L')[p^{\infty}]).
$$
Since $\mcal{C}$ is cofinal, 
every point involved in either 
$D_{A,K}(L)$ or $u(\hat{G}(L)[p^{\infty}])$
is defined over some $L'\in \mcal{C}$.
Taking the union over  $L'\in \mcal{C}$,
we obtain the desired result.
\end{proof}

\begin{proposition}
\label{pro-nilpotent}
Assume $e_K<p-1$. 
Let $L/K$ be a Galois extension with pro-nilpotent inertia.
Let $A_{/K}$ be an abelian variety with semistable reduction.
Then, condition $(C1)$ holds for $A_{/K}$ and $L/K$.
\end{proposition}

\begin{proof}
Let $\mcal{C}$ be the family of all finite Galois subextensions  of $L/K$.
Clearly $\mcal{C}$ is  cofinal.
Take any  $L'\in \mcal{C}$ and 
set $I:=I(L'/K)$.
If $\widehat{G}(L')[p]=0$,
then $z-1$ with $z=1\in Z(I)$ gives a bijection on $\widehat{G}(L')[p]$.
Suppose $\widehat{G}(L')[p]\not=0$.
Since $T$ splits over $K^{\mrm{ur}}$ and $e_K<p-1$,
we have $T(K^{\mrm{ur}})[p]=0$.
By \cite[Proposition 3.1]{ClXa08}, we also know that 
$\hat{B}(K^{\mrm{ur}})[p]=0$.
Hence  we have 
\begin{equation}\label{GL=0}
\widehat{G}(L')[p]^I=0
\end{equation}
by \eqref{rigid_ex2}.

Since $I$ is nilpotent, 
any $\ell$-Sylow subgroup $I_{\ell}$ of $I$ is normal.
Let $\Delta$ be the subgroup of $I$ generated by all $I_{\ell}$
for $\ell\not=p$.
Then we have a decomposition
$$
I=I_p\times \Delta.
$$
The natural projection from 
$I$ to the tame inertia group $I/I_p$ of $L'/K$
induces an isomorphism $\Delta\simeq I/I_p$.
In particular, $\Delta$ is cyclic,  and hence 
$$
\Delta\subset Z(I).
$$
Put $W=\widehat{G}(L')[p]$. 
Assume that $W^{\Delta}$ is not zero.
Then $W^{\Delta}$ is a non-zero $I_p$-stable submodule of $W$.
Since $I_p$ is a finite $p$-group, it follows that  
$W^I=(W^{\Delta})^{I_p}\not=0$.
This contradicts \eqref{GL=0}.
Hence $W^{\Delta}=0$.
Take any generator $z$ of the cyclic group $\Delta$.
We have $z\in \Delta\subset Z(I)$.
Since $W^{\Delta}=\ker(z-1\ \mrm{on}\ \widehat{G}(L')[p])$, $z-1$ is injective 
on $\widehat{G}(L')[p]$, and hence bijective.
Therefore, condition (C1) holds.
\end{proof}

Typical examples of Galois extensions $L/K$ with pro-nilpotent inertia are:
\begin{itemize}
\item[{\rm (i)}] abelian extensions;
\item[{\rm (ii)}] tamely ramified extensions;
\item[{\rm (iii)}] pro-$\ell$ Galois extensions for any prime $\ell$.
\end{itemize}

\begin{corollary}
\label{rigid_D}
Assume $e_K<p-1$. 
Let $L/K$ be a Galois extension with 
pro-nilpotent inertia.
Let $A_{/K}$ be an abelian variety.
\begin{itemize}
\item[{\rm (1)}] If $A$ has good reduction, then 
$D_{A,K}(L)$ coincides with the kernel of the reduction map on $A(L)[p^{\infty}]$.   
\item[{\rm (2)}] If $A$ is a $g$-dimensional abelian variety with split totally
degenerate reduction, then 
$D_{A,K}(L)$ is isomorphic to $\mu_{p^{\infty}}(L)^g$.
\end{itemize}
\end{corollary}
\begin{proof}
The result immediately follows from 
Example \ref{ex_D}, 
Proposition \ref{LD:semistable} (1) and Proposition \ref{pro-nilpotent}.
\end{proof}

\subsection{General cases}

In this section, we study local decomposition for abelian varieties under certain semistability restrictions. 
The goal of this section is to prove Theorem \ref{thm:intro-local}.


First we prove the following descent theorem.
It is needed for the proof of Theorem \ref{thm:intro-local},
but is also of independent interest and plays an important role
in our global applications, such as Theorem~\ref{cor:cyclotomic-rigidity_intro}.
    
\begin{theorem}
\label{Yoshida_descent}
Let $L/K$ be a totally ramified pro-$p$
Galois extension.
Let $A_{/K}$ be an abelian variety
which acquires semistable reduction over a finite extension $K'/K$ satisfying $e_{K'}<p-1$.
Then,
$$
A(L)_{\mrm{tor}}=A(K)_{\mrm{tor}}.
$$
\end{theorem}

\begin{lemma}\label{Yoshida_descent_lemma}
Let $L/K$ be a totally ramified pro-$p$ Galois extension.
Let $K'/K$ be a finite extension whose ramification index $e(K'/K)$
is prime to $p$.
Put $L':=LK'$. 
Then,
$L'/K'$ is a totally ramified pro-$p$ Galois extension and $L\cap K'=K$.
\end{lemma}
\begin{proof}
Let $M/K$ be a finite Galois subextension in $L/K$
and put $M'=MK'$.
Since $L/K$ is totally ramified pro-$p$, 
we have $[M:K]=e(M/K)=p^n$ for some $n\ge 0$.
Since $e(M/K)$ divides $e(M'/K)=e(M'/K')e(K'/K)$
and $e(K'/K)$ is prime-to-$p$,
$p^n$ divides $e(M'/K')$.
Thus 
$$
p^n\le e(M'/K')\le [M':K']=[M:M\cap K']\le [M:K]=p^n.
$$
Therefore all the inequalities are equalities. 
In particular,
$M\cap K'=K$,
which implies $L\cap K'=K$.
Furthermore, $e(M'/K')=[M':K']$,
which implies that $L'/K'$ is totally ramified and pro-$p$.
\end{proof}

\begin{proof}[{Proof of Theorem \ref{Yoshida_descent}}] 
We divide the proof into three steps.

\medskip 

\noindent
{\bf Step 1.} Assume that $A/K$ has semistable reduction and that $e_K<p-1$.
We show 
$$
A(L)[p^{\infty}]=A(K)[p^{\infty}].
$$
Since $L/K$ is totally ramified and pro-$p$, its inertia group is
pro-nilpotent. Hence condition (C1) is satisfied.
By Proposition \ref{LD:semistable}, 
we have 
$A(L)[p^{\infty}]=D_{A,K}(L)\oplus A(K)[p^{\infty}]$.
Taking $\mrm{Gal}(L/K)$-fixed parts, 
we see that $D_{A,K}(L)^{\mrm{Gal}(L/K)}=0$.
We claim that this implies $D_{A,K}(L)=0$.
Indeed, if $D_{A,K}(L)\ne0$, then
$D_{A,K}(L)[p]\ne0$.
The action of the pro-$p$ group $\Gal(L/K)$ on the finite
$\mathbb F_p$-vector space $D_{A,K}(L)[p]$ factors through a finite
$p$-group, and hence
$D_{A,K}(L)[p]^{\Gal(L/K)}\ne0$,
a contradiction.
Thus $D_{A,K}(L)=0$, proving
\[
A(L)[p^\infty]=A(K)[p^\infty]
\]
as desired.

\medskip 

\noindent
{\bf Step 2.} Assume that $A/K$ has semistable reduction.
We show that, for every prime $\ell\not=p$,
$$
A(L)[\ell^{\infty}]=A(K)[\ell^{\infty}].
$$
Since $A_{/K}$ has semistable reduction,
$I_K$ acts unipotently on $A(L)[\ell^{\infty}]$,
and hence this $I_K$ action factors through a pro-$\ell$ quotient.
On the other hand, since $L/K$ is totally ramified, the action of
$I_K$ on $A(L)[\ell^\infty]$ factors through
$\mrm{Gal}(L/K)$,
which is pro-$p$.
Hence we obtain that $\mrm{Gal}(L/K)$ acts trivially on 
$A(L)[\ell^{\infty}]$,
which gives the desired result.

\medskip 

\noindent
{\bf Step 3.} We prove the theorem in general.
Let $K'/K$ be as in the statement
and put $L'=LK'$. 
By Lemma \ref{Yoshida_descent_lemma},
$L'/K'$ is a totally ramified pro-$p$ extension
and $L\cap K'=K$.
Since $A$ has semistable reduction over $K'$ and 
$e_{K'}<p-1$, applying Step 1 and Step 2 over $K'$ gives 
$A(L')_{\mrm{tor}}=A(K')_{\mrm{tor}}$.
Therefore, 
$$
A(L)_{\mrm{tor}}=A(L)_{\mrm{tor}}\cap A(L')_{\mrm{tor}}
=A(L)_{\mrm{tor}}\cap A(K')_{\mrm{tor}}
=A(L\cap K')_{\mrm{tor}}
=A(K)_{\mrm{tor}}.
$$
This proves the theorem.
\end{proof}

Next we prove Theorem \ref{thm:intro-local}.

\begin{lemma}\label{LD:general_lemma}
Let the notation be as in Theorem \ref{thm:intro-local}.
For any finite Galois subextension $L'/K$ of $L/K$ and any $n>0$,
put 
$$
D_{A,K}(L')_n=\sum_{\sigma\in I(L'/K)} (\sigma-1)A(L')[p^n].
$$
Then, for any $I(L'/K)$-stable submodule $N'$ of $A(L')[p^n]$,
we have 
$$
N'=(N'\cap D_{A,K}(L')_n)\oplus N'^{I(L'/K)}.
$$
\end{lemma}

\begin{proof}
Set $I:=I(L'/K)$.
Since $I$ is nilpotent, we have a decomposition
$I=I_p\times \Delta$
where $I_p$ is a finite $p$-group and 
$\Delta$ is a group of order prime-to-$p$.
Note that $\Delta$ is cyclic since it is 
isomorphic to the tame inertia quotient of  $I$.
In particular, $\Delta\subset Z(I)$.

Set $M:=A(L')[p^n]$ and $L'_0:=(L')^I$.
Then $L'_0/K$ is unramified, and
$L'^\Delta/L'_0$ is a totally ramified Galois extension of
$p$-power degree.
By the assumption of Theorem~\ref{thm:intro-local}, there is a finite
extension $K'/K$ over which $A$ has semistable reduction and such that
$e_{K'}<p-1$.
Note that $A$ has semistable reduction over
$K'L'_0$ and $e_{K'L'_0}=e_{K'}<p-1$.
Applying Theorem \ref{Yoshida_descent}
to $A_{/L'_0}$, we have 
\begin{equation}
\label{delta=I}
M^{\Delta}=M^I.
\end{equation}
Now we define 
$$
e_{\Delta}\colon M\to M
$$
by $e_{\Delta}(x)=|\Delta|^{-1}\sum_{\delta\in \Delta} \delta(x)$.
This is well defined since $p\nmid |\Delta|$.
By definition, we see
\begin{equation}
\label{e_delta1}
e_{\Delta}=\mrm{id} \ \mrm{on}\ M^{\Delta}\quad \mrm{and}\quad 
e_{\Delta}(M)=M^{\Delta}.
\end{equation}
Furthermore, we have 
\begin{equation}
\label{e_delta2}
\ker e_{\Delta}=D_{A,K}(L')_n.
\end{equation}
Indeed, let $x\in M$. Since $\Delta\subset Z(I)$,
we have $e_{\Delta}(\sigma-1)(x)=(\sigma-1)(e_{\Delta}(x))=0$.
On the other hand,
if $x\in \ker e_{\Delta}$, 
$$
x=x-e_{\Delta}(x)
=-|\Delta|^{-1}\sum_{\delta\in \Delta} (\delta-1)(x)\in D_{A,K}(L')_n.
$$
Thus \eqref{e_delta2} follows.

Since $N'$ is $I$-stable, we have $e_{\Delta}(N')\subset N'$. 
Hence we obtain 
$$
N'=\ker (e_{\Delta} \mrm{on}\ N')\oplus e_{\Delta}(N').
$$
Now we have 
$\ker (e_{\Delta} \mrm{on}\ N')=N'\cap \ker e_{\Delta}=N'\cap D_{A,K}(L')_n$
by \eqref{e_delta2}.
Furthermore, we also have  $e_{\Delta}(N')=N'^I$
by \eqref{delta=I}.
This completes the proof.
\end{proof}

\begin{proof}[{Proof of Theorem 
\ref{thm:intro-local}}]
For the proof of Theorem 
\ref{thm:intro-local}, it suffices to show the following.
\begin{itemize}
\item[{\rm (I)}] 
For every $I(L/K)$-stable submodule
$N\subset A(L)[p^\infty]$,
\[
 N=\bigl(N\cap D_{A,K}(L)\bigr)\oplus N^{I(L/K)}.
\]
\item[{\rm (II)}] We have 
$A(L)_{p'}=A(L\cap K^{\mrm{tame}})_{p'}$.
\end{itemize}

\medskip

First we show (I). For any $n>0$, put 
$$
N_n=N[p^n]=N\cap A(L)[p^n]. 
$$
Since $N_n$ is finite, we may choose 
a finite Galois subextension $L_n/K$ of $L/K$ 
such that $N_n\subset A(L_n)$.  
The action of $I(L/K)$ on $N_n$ factors through $I(L_n/K)$,
and $N_n$ is an $I(L_n/K)$-stable submodule of $A(L_n)[p^n]$.
By Lemma \ref{LD:general_lemma}, we have 
\begin{equation}
\label{Nn2}
N_n=(N_n\cap D_{A,K}(L_n)_n)\oplus N_n^{I(L_n/K)}.
\end{equation}
By the same argument as in the proof of
Proposition~\ref{LD:semistable}(2),
with $u(\widehat{G}(L_n)[p^n])$ replaced by
$D_{A,K}(L_n)_n$, equation \eqref{Nn2} yields 
$N=(N\cap D_{A,K}(L))\oplus N^{I(L/K)}$.

\medskip

Next we show (II). Take  any prime $\ell\not=p$.
By assumption, $A$ acquires semistable reduction 
over a tamely ramified finite extension $K'/K$.
Since $A$ has semistable reduction over $K'$,
$I_{K'}$ acts unipotently on $A(\overline{K})[\ell^{\infty}]$,
and hence this $I_{K'}$-action factors through a pro-$\ell$ quotient.
Hence the wild inertia $P_{K'}$ acts trivially on 
$A(\overline{K})[\ell^{\infty}]$.
On the other hand, we have $G_{K'}\supset P_K$ since $K'/K$ is tamely ramified.
Thus $P_{K'}=P_K$, and therefore $P_K$ acts trivially on $A(\overline{K})[\ell^{\infty}]$.
As this holds for every $\ell\ne p$, $P_K$ acts trivially on
$A(\overline K)_{p'}$.
Consequently, 
$A(L)_{p'}=A(L\cap K^{\mrm{tame}})_{p'}$
as desired.
\end{proof}

\subsection{The case of good reduction}
\label{secton:LD_CDVF}

In this section, we study local decomposition 
for an abelian variety $A_{/K}$ with good reduction. 
In this case, local decomposition 
behaves well.

We denote by $\widehat{A}$ and $\bar{A}$ the formal group over 
$\mcal{O}_K$ associated with $A$ and 
the reduction of $A$, respectively.
Throughout this paper, for an algebraic extension $L$ of 
$K$ with residue field $k_L$ and maximal ideal $\mbf{m}_L$,
we denote by $\widehat{A}(L)$ the group $\widehat{A}(\mbf{m}_L)=\mbf{m}_L^{h}$ 
where the group structure on $\widehat{A}(L)$ is 
induced by the formal group law  of $\widehat{A}$. Here, $h$ is the dimension of $\widehat{A}$.
Let $L/K$ be a Galois extension
with residue field $k_L$. 
By \cite[Part C, Appendix C2, Theorem 2.6]{HiSi00},
we may identify $\widehat{A}(L)$ with the kernel of the reduction map 
$A(L)\to \bar{A}(k_L)$; that is,
$$
\widehat{A}(L)=\mrm{ker}\left(A(L)\to \bar{A}(k_L)\right).
$$
In the next section, we study the formal part
$\widehat{A}(L)[p^{\infty}]$
under further assumptions.

\medskip

The condition $(C1)$ 
for $A_{/K}$ and $L/K$ is now described as follows.

\begin{itemize}
\item[{\bf (C1):}] 
There exists a cofinal family $\mathcal{C}$ of finite Galois subextensions
$L'/K$ of $L/K$ such that, for every $L'\in\mathcal{C}$, 
one can choose $z\in Z(I(L'/K))$
so that 
$z-1$ acts bijectively on $\widehat{A}(L')[p]$.
\end{itemize}

In the case of good reduction, the local decomposition theorem takes the following simple form.

\begin{theorem}\label{LD:p-adic}
Let $K$ be a complete discretely valued field of mixed characteristic
$(0,p)$ with perfect residue field.
Let $L/K$ be a Galois extension.
Let $A$ be an abelian variety over $K$ with good reduction.
Assume that the condition $(C1)$ holds for $A_{/K}$ and $L/K$.
Then, 
\[
A(L)_{\mrm{tor}}
= 
\widehat{A}(L)[p^{\infty}]
\oplus A(L\cap K^{\mrm{ur}})_{\mrm{tor}}.
\]
Moreover, for every $I(L/K)$-stable submodule
$N\subset A(L)[p^\infty]$,
\[
 N=\bigl(N\cap \widehat{A}(L)[p^{\infty}]\bigr)\oplus N^{I(L/K)}.
\]
\end{theorem}

\begin{proof}
The assertions concerning the $p$-primary parts follow from
Proposition~\ref{LD:semistable}, since in the good reduction case, we have
$
u(\hat{G}(L)[p^{\infty}])=\widehat{A}(L)[p^\infty].
$
On the other hand, the N\'eron-Ogg-Shafarevich criterion gives
$A(L)_{p'}=A(L\cap K^{\mrm{ur}})_{p'}$.
Combining the $p$-primary and prime-to-$p$ parts gives the first
assertion, while the second assertion follows directly from
Proposition~\ref{LD:semistable}.
\end{proof}

\section{Formal torsion over totally ramified extensions}
Let $K$ be a complete discretely valued field of 
mixed characteristic $(0,p)$ with perfect residue field $k$.
For any finite extension $K'/K$, 
we denote by $e_{K'}$ the absolute ramification index of $K'$
and also denote by $e^{\mrm{tame }}_{K'}$ the tame ramification index of 
$K'/K$, that is, the prime-to-$p$ part of 
the relative ramification index of $K'/K$.
Let $g>0$ be an integer.
Let $A_{/K}$ be a $g$-dimensional abelian variety 
with {\it semistable reduction}.
Let 
$$
0\to T\to G\overset{v}{\to} B\to 0 \quad \mrm{and}\quad 
0\to N^{\mrm{an}}\to G^{\mrm{an}}\overset{u}{\to} A^{\mrm{an}}\to 0
$$
be the non-Archimedian rigid uniformization
as in Section \ref{LD_st}. 
We have seen that $u(\widehat{G}(L)[p^{\infty}])$
plays an important role in the local decomposition.
Our aim in this section is to give an explicit description of the formal part 
$\widehat{G}(L)[p^{\infty}]$.

We need some notation.
Let
$$
0\to \mcal{T}\to \mcal{G}\to \mcal{B}\to 0
$$
be the Raynaud extension over $\cO_K$, where $\mcal{T}$ is a torus,
$\mcal{B}$ is an abelian scheme, and $\mcal{G}$ is a semi-abelian scheme,
whose generic fiber is
$
0\to T\to G\overset{v}{\longrightarrow} B\to 0
$
above.
This induces an exact sequence of finite flat $\cO_K$-group schemes 
$$
0\to \mcal{T}[p^n]\to \mcal{G}[p^n]\overset{v}{\to} \mcal{B}[p^n]\to 0.
$$
For any finite flat group scheme $H$ over $\cO_K$, we denote by
$H^0$ its identity component and by $H^{\aet}$ its maximal \'etale quotient.
Furthermore, we denote by $H^D$ the Cartier dual of $H$.
In what follows, exact sequences of finite flat
group schemes over $\cO_K$ are understood to be
exact as sequences of fppf sheaves.

\subsection{Formal torsion via finite flat group schemes}
\label{Sec_FG_ur}

We study the formal part 
$\widehat{G}(L)[p^{\infty}]$.

\subsubsection{The condition (C2) and its sufficiency}

We consider the following condition (C2) for $A_{/K}$ and $L/K$.

\begin{itemize}
\item[{\bf (C2):}]
For every $n>0$, the natural closed immersion
$\left(\mcal{B}^{\vee}[p^n]^{\aet}\right)^D
\hookrightarrow
\mcal{B}[p^n]$
induces a $G_K$-equivariant identification
$$
\widehat{B}(L)[p^n]
=
\left(\mcal{B}^{\vee}[p^n]^{\aet}\right)^D(L),
$$
compatible with the natural transition maps in $n$.
\end{itemize}

\noindent
The abelian scheme $\mcal{B}$ is canonically
determined by $A$ up to isomorphism
(see \cite[Expos\'e IX, Section 7.2]{SGA7}).
Hence the condition above is independent
of the choice of the Raynaud uniformization.
Note that, if (C2) is satisfied, then we have  a natural isomorphism
$$
\left(\mcal{B}^{\vee}[p^n]^{\aet}\right)^D(L)
\simeq 
\mrm{Hom}_{G_L}(\bar{B}^{\vee}(\bar k)[p^n], \mu_{p^n}(\overline{K}))
$$
of $G_K$-modules by definition of Cartier duals.
Furthermore, 
there exists a $G_K$-equivariant isomorphism
$$
\widehat{B}(L)[p^\infty]
\simeq \mrm{Hom}_{\mrm{cont}, G_L}(T_p(\bar{B}^{\vee}), \mu_{p^{\infty}}(\overline{K})).
$$
Indeed, we have natural $G_K$-equivariant isomorphisms
$$
T_p(\bar B^\vee)/p^nT_p(\bar B^\vee)
\simeq
\bar B^\vee(\bar k)[p^n].
$$
Moreover, every continuous homomorphism from
$T_p(\bar B^\vee)$ to the discrete group $\mu_{p^\infty}$
has finite image and hence factors through
$T_p(\bar B^\vee)/p^nT_p(\bar B^\vee)$ for some $n$.
Thus, taking the union over $n$ in the isomorphisms of (C2)
gives the desired isomorphism.

\begin{proposition}
\label{C2}
Let $L/K$ be a Galois extension and  
let $A$ be an abelian variety over $K$ with good reduction.
Put $g=\dim A$.
Assume that one of the following holds:
\begin{itemize}
\item[{\rm (a)}] $e_K=1$, 
$L/K$ is totally ramified,  and the tame ramification index of $L/K$ divides $p-1$.
\item[{\rm (b)}] $L/K$ is totally ramified, and $K/\mbb{Q}_p$ is a finite unramified 
extension of degree $f$ with
$\mrm{gcd}(f,(2g)!)=1$.
\item[{\rm (c)}] $A$ has ordinary reduction.
\end{itemize}
Then condition {\rm (C2)} holds  for $A_{/K}$ and $L/K$.
\end{proposition}

\begin{proof}
Fix an integer $n>0$ and 
put $\mcal{H}=\mcal{A}[p^n]$.
Then $\mcal{H}^D\simeq \mcal{A}^{\vee}[p^n]$, 
where $\mcal{A}^{\vee}$ is the dual abelian scheme of $\mcal{A}$.
We also put 
$$
\mcal{M}=(\mcal{H}^D)^{\aet}=(\mcal{A}^{\vee}[p^n])^{\aet}.
$$
Dualizing the natural surjection 
$\mcal{H}^D\to \mcal{M}$, we may regard 
$\mcal{M}^D$ as a closed finite flat subgroup scheme of $\mcal{H}$.
Since there is a canonical isomorphism $\mcal{H}^0(L)\simeq \widehat{A}(L)[p^n]$,
it suffices to prove 
$$
\mcal{M}^D(L)=\mcal{H}^0(L)
$$
as subgroups of $A(\overline{K})[p^n]$.

\medskip

We first show that
$\mcal{M}^D$ is a closed subgroup scheme
of $\mcal{H}^0$.
Taking the Cartier dual of the connected-\'etale sequence of $\mcal{H}$ gives 
an exact sequence of finite flat $\cO_K$-group schemes 
$$
0\to (\mcal{H}^{\aet})^D\to \mcal{H}^D\to (\mcal{H}^0)^D\to 0.
$$
Since $(\mcal{H}^{\aet})^D$ is multiplicative and is of $p$-power order,
it is connected, and  the maximal \'etale quotient $\mcal{M}$ of
$\mcal{H}^D$ is naturally isomorphic to that of $(\mcal{H}^0)^D$. 
Thus, the natural surjection $\mcal{H}^D\to \mcal{M}$ factors through $(\mcal{H}^0)^D$.
By Cartier duality, the closed immersion 
$\mcal{M}^D\to \mcal{H}$ factors through  $\mcal{H}^0$.
Hence $\mcal{M}^D$ is a closed subgroup scheme of  $\mcal{H}^0$.

\medskip

We first consider case (c).
Since the dual abelian variety 
$A^{\vee}$ has good ordinary reduction, we have 
$$
\mrm{rank}(\mcal{M}^D)=\mrm{rank}(\mcal{H}^0)=p^{ng}.
$$
Hence 
the closed immersion 
$\mcal{M}^D\hookrightarrow \mcal{H}^0$ of finite flat group schemes 
is in fact an isomorphism
\begin{equation}\label{ordisom}
\mcal{M}^D\simeq \mcal{H}^0.
\end{equation}
This proves (c).

\medskip 

Next we consider case (a).
Put $\mcal{Q}:=\mcal{H}^0/\mcal{M}^D$.
It suffices to show that $\mcal{Q}(L)=0$.
For this, first we note that $\mcal{Q}$ and $\mcal{Q}^D$ are connected.
Indeed, since $\mcal{H}^0$ is connected, so is $\mcal{Q}$.
Moreover, by the above identification of $\mcal{M}$
with the maximal \'etale quotient of $(\mcal{H}^0)^D$,
we have
$$
\mcal{Q}^D
\simeq
\ker((\mcal{H}^0)^D\to \mcal{M})
\simeq
\ker((\mcal{H}^0)^D\to ((\mcal{H}^0)^D)^{\aet})
=
((\mcal{H}^0)^D)^0.
$$
Thus $\mcal{Q}^D$ is also connected.

Assume that
$\mcal{Q}(L)\ne0$, and put $V=\mcal{Q}(L)[p]\ne0$.
By assumption on $L/K$, we have an exact sequence 
$$
1\to P\to \mrm{Gal}(L/K)\to \Delta\to 1
$$
where $P$ is the wild inertia of $\mrm{Gal}(L/K)$ and
$\Delta$ is the tame inertia quotient of $\mrm{Gal}(L/K)$.
Note that $\Delta$ is cyclic.
Since $P$ is pro-$p$, we have $V^P\not=0$.
As $\Delta$ is cyclic of order dividing $p-1$, there is a
$\Delta$-stable one-dimensional $\mbb F_p$-subspace
$
W\subset V^P.
$
Thus $W$ is $G_K$-stable.
Let $\mcal{C}_K\subset \mcal{Q}_K$ be the finite \'etale subgroup scheme corresponding to
$W$, and let $\mcal{C}$ be the scheme-theoretic closure of $\mcal{C}_K$ in $\mcal{Q}$.
Then $\mcal{C}$ is a finite flat group scheme of order $p$.
Moreover, $\mcal{C}$ is connected,
since its special fiber is a closed
subgroup scheme of the connected finite
group scheme $\mcal{Q}_k$.

Since $e_K=1$, it follows from the classification of Oort-Tate \cite[Theorem 2 and Remarks 5.6]{TaOo70} that 
$\mcal{C}$ is of multiplicative type. Thus $\mcal{C}^D$ is \'etale.
On the other hand, the closed immersion  $\mcal{C}\hookrightarrow \mcal{Q}$ induces a
surjection
$
\mcal{Q}^D\to \mcal{C}^D,
$
which is impossible since $\mcal{Q}^D$ is connected
 and $\mcal{C}^D$ is a non-trivial \'etale group scheme.
Therefore $\mcal{Q}(L)=0$.
This proves (a).

\medskip 

Finally, we consider case (b).
Put $N=\widehat{A}(L)[p^n]$.
Let $L_n (\subset L)$ be the field of definition $K(N)$ of the points of $N$. 
Then, $L_n/K$ is a finite totally ramified Galois extension, 
and $\mrm{Gal}(L_n/K)$ acts faithfully on $N$. 
Setting $d_n:=\dim_{\mbb{F}_p} N/pN$,
we have $d_n\le 2g$. Thus $\mrm{gcd}(f,d_n!)=1$ by assumption (b).
Applying Lemma \ref{YoshidasLemma} below,
the tame ramification index of $L_n/K$
divides $p-1$.
Therefore, by case (a),
we obtain
$$
N=\widehat{A}(L_n)[p^n]= \mcal{M}^D(L_n)\subset \mcal{M}^D(L). 
$$
Since $\mcal{M}^D$ is a closed
subgroup scheme of $\mcal{A}[p^n]^0$, 
we obtain $N= \mcal{M}^D(L)$.
This completes the proof in the case (b).
\end{proof}

In the proof above, 
we used the following lemma.
\begin{lemma}
\label{YoshidasLemma}
Let $K/\mbb{Q}_p$ be the unramified extension of degree $f$, and 
let $L/K$ be a finite totally ramified Galois
extension.
Assume that $\mrm{Gal}(L/K)$ acts faithfully on a finite module $N$ of $p$-power order.
Set $d:=\mrm{dim}_{\mbb{F}_p} N/pN$.
If $\mrm{gcd}(f, d!)=1$,
then the tame ramification index of $L/K$
divides $p-1$.
\end{lemma}

\begin{proof}
Let $P$ be the wild inertia subgroup of $\mrm{Gal}(L/K)$,
and denote by $m$ the tame ramification index of $L/K$. 
Let $\varpi$ be a uniformizer of $L$
and denote by $\mbf{m}_L$ the maximal ideal of $\cO_L$.
Since $L/K$ is totally ramified,
the map 
$$
\mrm{Gal}(L/K)\to \mbb{F}_{p^f}^{\times}
$$
given by $\sigma\mapsto \sigma(\varpi)/\varpi$ mod $\mbf{m}_L$ is well-defined.
Furthermore, this is a homomorphism with kernel $P$. Thus $m$ divides $p^f-1$.

The kernel of the natural homomorphism 
$\iota:\mrm{Aut}(N)\to \mrm{Aut}(N/pN)
=GL_{\mbb{F}_p}(N/pN)$ is 
of $p$-power order.
Indeed, suppose $p^{\ell}N=0$, 
and let 
$$
U_i:=\{\sigma\in \mrm{Aut} (N)\mid (\sigma-1)N\subset p^iN\}.
$$
Then, $U_1=\ker \iota$, $U_{\ell}=1$
and the homomorphism 
$$
U_i/U_{i+1}\to \mrm{Hom}(N,p^iN/p^{i+1}N)
$$
defined by $\sigma\mapsto (x\mapsto (\sigma-1)x)$ is injective.
Since the right-hand side is a finite abelian $p$-group,
each quotient $U_i/U_{i+1}$ has $p$-power order.
Thus $U_1$ is of $p$-power order as desired.

Take a lift $\sigma$ of a generator of $\mrm{Gal}(L/K)/P$. Since the order of $\sigma$ is of the form $mp^a$, the order of $\delta:=\sigma^{p^a}$
is $m$. Since $m$ is prime-to-$p$, $\mrm{Gal}(L/K)$ acts faithfully on $N$ and the kernel of $\iota$ is of $p$-power order, 
it follows that 
the image $\bar{\delta}$ of  $\delta$ 
in $GL_{\mbb{F}_p}(N/pN)$ is of order $m$.
Thus the minimal polynomial $f(T)$ of $\bar{\delta}$
divides $T^m-1$. Since $m$ is prime to $p$,
$\bar{\delta}$ is diagonalizable (over $\bar{\mbb{F}}_p$).
Since $m$ divides $p^f-1$,
all the eigenvalues $\lambda$ of $\bar{\delta}$
are contained in $\mbb{F}_{p^f}$.
Thus, the degree of the irreducible factor 
of $f(T)$ corresponding to each 
eigenvalue $\lambda$ satisfies 
$$
d_{\lambda}\mid f\quad \mrm{and}\quad 
1\le d_{\lambda}\le d.
$$
By $\mrm{gcd}(f,d!)=1$, we have $d_{\lambda}=1$.
This implies that all the eigenvalues of
$\bar{\delta}$ belong to $\mbb{F}_p^\times$.
Since $\bar{\delta}$ is diagonalizable,
we have $\bar{\delta}^{p-1}=1$.
As $\bar{\delta}$ has order $m$,
we conclude that $m\mid p-1$.
This completes the proof.
\end{proof}

\begin{corollary}
\label{ordrep}
Assume that  $A_{/K}$ has good ordinary reduction. 
Then, there exists a natural isomorphism 
$$
T_p(\widehat{A})\simeq T_p(\bar{A}^{\vee})^{\vee}(1)
$$
of $G_K$-modules.
\end{corollary}

\begin{proof}
The result follows from \eqref{ordisom}
and Cartier duality.
\end{proof}

\subsubsection{Formal torsion and Galois actions}

\begin{theorem}
\label{hatG}
Let $L/K$ be a Galois extension. 
Let $A$ be an abelian variety over $K$ with semistable reduction.
Assume that {\rm (C2)} holds for $A_{/K}$ and $L/K$.
Then, for every $n>0$, there is a natural $G_K$-equivariant
isomorphism
$$
\widehat G(L)[p^n]
\simeq
\left((\mcal G[p^n]^D)^{\aet}\right)^D(L),
$$
compatible with the natural transition maps in $n$.
Furthermore, there is a natural exact sequence of finite flat
$\cO_K$-group schemes
$$
0\to \mcal T[p^n]
\to
\left((\mcal G[p^n]^D)^{\aet}\right)^D
\to
\left((\mcal B^\vee[p^n])^{\aet}\right)^D
\to0,
$$
compatible with the natural transition maps in $n$.\end{theorem}
\begin{proof}
Set $\mcal{M}_{G,n}:=(\mcal{G}[p^n]^D)^{\aet}$ and 
$\mcal{M}_{B,n}:=(\mcal{B}[p^n]^D)^{\aet}$.
Then $\mcal{M}_{G,n}^D$ and $\mcal{M}_{B,n}^D$
are closed finite flat subgroup schemes of  $\mcal{G}[p^n]$ and $\mcal{B}[p^n]$,
respectively.
We also set $\mcal{H}_n:=\mcal{G}[p^n]\times_{\mcal{B}[p^n]} \mcal{M}_{B,n}^D$.
Then $\mcal H_n$ is a closed finite flat subgroup scheme of
$\mcal G[p^n]$, and by definition we have a natural exact sequence
\begin{equation}\label{THM}
0\to\mcal T[p^n]
\to\mcal H_n
\to\mcal M_{B,n}^D
\to0.
\end{equation}

We claim 
\begin{equation}\label{MH}
\mcal{M}_{G,n}^D=\mcal{H}_n.
\end{equation}
Dualizing~\eqref{THM}, we obtain
$0\to \mcal M_{B,n}
\to \mcal H_n^D
\to \mcal T[p^n]^D
\to 0.$
Since both $\mcal{M}_{B,n}$ and $\mcal{T}[p^n]^D$ are \'etale, 
$\mcal{H}_n^D$ is also \'etale.
Hence the dual of 
the closed immersion $\mcal{H}_n\hookrightarrow \mcal{G}[p^n]$
factors through the maximal \'etale quotient $\mcal{M}_{G,n}$ of $\mcal{G}[p^n]^D$.
Dualizing again, we obtain
$\mcal H_n\subset\mcal M_{G,n}^D$. 
Conversely,
since $\mcal{M}_{G,n}$ is \'etale,
the composite $\mcal{B}[p^n]^D\hookrightarrow \mcal{G}[p^n]^D\to \mcal{M}_{G,n}$ 
factors through the maximal \'etale quotient $\mcal{M}_{B,n}$.
Dualizing this factorization shows that the closed immersion
$\mcal M_{G,n}^D\hookrightarrow\mcal G[p^n]$
factors through $\mcal H_n$.
Thus, $\mcal{M}_{G,n}^D\subset \mcal{H}_n$,
which proves~\eqref{MH}.

It remains to prove $\widehat{G}(L)[p^n]=\mcal{M}_{G,n}^D(L)$.
For any $x\in \hat{G}(L)[p^n]$,
we have $v(x)\in \widehat{B}(L)[p^n]$,
which gives $v(x)\in \mcal{M}_{B,n}^D(L)$ by (C2).
This implies $x\in \mcal{H}_{n}(L)$.
Hence it follows from the claim above that 
$\hat{G}(L)[p^n]\subset \mcal{H}_n(L)=\mcal{M}_{G,n}^D(L)$.
Conversely, take any $x\in \mcal{M}_{G,n}^D(L)$.
By the claim, we have $x\in \mcal{H}_n(L)$.
By (C2), $v(x)\in \mcal{M}_{B,n}^D(L)=\widehat{B}(L)[p^n]$.
This implies $x\in \widehat{G}(L)[p^n]$. Thus we obtain
$\hat{G}(L)[p^n]\supset \mcal{H}_n(L)=\mcal{M}_{G,n}^D(L)$.
This proves $\widehat{G}(L)[p^n]=\mcal{M}_{G,n}^D(L)$ as desired.

The equality
\[
\widehat G(L)[p^n]=\mcal M_{G,n}^D(L)
\]
is an equality of $G_K$-submodules of
$\mcal G[p^n](\overline K)$.
Moreover, the constructions of
$\mcal M_{G,n}$, $\mcal M_{B,n}$, and $\mcal H_n$
are functorial with respect to the natural transition maps in $n$.
Together with the compatibility assumed in {\rm (C2)}, this shows that
these identifications are compatible with the natural transition maps in $n$.

Finally, the last assertion follows from
\eqref{THM} and \eqref{MH}.
Its compatibility with the natural transition maps in $n$
follows from the same functoriality.
\end{proof}

\begin{corollary}\label{hatG_cyclic}
Let $A_{/K}$ and $L/K$ be as in Theorem \ref{hatG}.
Then the following assertions hold.
\begin{itemize}
\item[{\rm (1)}]
The group $I_K$ acts on $\widehat{G}(L)[p^\infty]$
via the cyclotomic character $\chi_p$.
Furthermore, if $L/K$ is totally ramified,
then the action of $G_K$ on
$\widehat{G}(L)[p^\infty]$ is scalar.
\item[{\rm (2)}]
If $L/K$ is totally ramified and
$L$ contains $K(\mu_{p^\infty})$,
then $G_K$ acts on
$\widehat{G}(L)[p^\infty]$ via $\chi_p$.
\end{itemize}
\end{corollary}
\begin{proof}
By Theorem \ref{hatG} and Cartier duality, it is enough to consider the natural
$G_K$-action on
$$
\mrm{Hom}_{G_L}((\mcal{G}[p^n]^D)^{\aet}(\bar k), \mu_{p^n}(\overline{K})).
$$
Since $I_K$ acts trivially on $(\mcal{G}[p^n]^D)^{\aet}(\bar k)$, 
the action of $I_K$ on the above Hom group
comes entirely from its action on $\mu_{p^n}(\overline{K})$, which is given
by $\chi_p$.
Now assume that $L/K$ is totally ramified.
Take any $\sigma\in G_K$
and any $P\in\widehat{G}(L)[p^\infty]$.
Since $L/K$ is totally ramified,
there exists $\tau\in I_K$ such that
$\sigma|_L=\tau|_L$.
Since $P$ is $L$-rational, we have  $\sigma(P)=\tau(P)=\chi_p(\tau)P$.
This proves (1).
Assume moreover $K(\mu_{p^\infty})\subset L$. 
Then $\sigma|_L=\tau|_L$ implies
$\sigma^{-1}\tau\in G_L$.
Since $\chi_p$ is trivial on $G_L$, we have
$\chi_p(\sigma)=\chi_p(\tau)$, which proves (2).
\end{proof}

Passing to the maximal unramified extension gives the following useful consequence.
In Section \ref{sec:applications}, we apply this result to the study of torsion over the maximal abelian extension, including explicit determinations for Jacobians.

\begin{corollary}\label{hatG_ab}
Assume that $e_K=1$.
Let $A_{/K}$ be an abelian variety with semistable reduction.
Then, the group 
$I_K$ acts on $\widehat{G}(K^{\mrm{ur}}(\mu_{p^{\infty}}))[p^{\infty}]$ via $\chi_p$.
\end{corollary}

\begin{proof}
For a finite unramified extension $K'/K$,
set $L':=K'(\mu_{p^{\infty}})$.
Then $L'/K'$ is totally ramified and its tame ramification index is $p-1$.
By Proposition \ref{C2},
(C2) holds for $A_{/K'}$ and $L'/K'$.
Thus Corollary \ref{hatG_cyclic} implies that $I_K=I_{K'}$
acts on $\widehat{G}(L')[p^{\infty}]$ via $\chi_p$.
Varying $K'$, we obtain the desired result.
\end{proof}

\begin{corollary}\label{hatG_p-adic}
Let $K$ be a $p$-adic field, 
let $L/K$ be a totally ramified Galois extension,
and let $A_{/K}$ have semistable reduction.
Assume (C2) and $L\supset K(\mu_{p^{\infty}})$. 
Then, there exists a non-canonical isomorphism of groups 
$$
\widehat{G}(L)[p^{\infty}]\simeq 
\varinjlim_n \ (\mcal{G}[p^n]^D)^{\aet}(k).
$$
In particular, if $A$ has good reduction,
there exists a non-canonical isomorphism of groups
$$
\widehat{A}(L)[p^{\infty}]\simeq 
\bar A^{\vee}(k)[p^{\infty}].
$$
\end{corollary}

\begin{proof}
Set $\mcal{M}_{G,n}:=(\mcal{G}[p^n]^D)^{\aet}$.
Let $\mcal{M}_G$ be the \'etale $p$-divisible group defined by 
$(\mcal{M}_{G,n})_n$ with the natural 
transition maps $\mcal{M}_{G,n}\hookrightarrow \mcal{M}_{G,n+1}$.
We denote by $T_p:=T_p(\mcal{M}_{G})$ the Tate module of $\mcal{M}_G$,
and set $V_p:=T_p\otimes_{\mbb{Z}_p}\mbb{Q}_p$.
By Theorem \ref{hatG} and Cartier duality, 
we have natural isomorphisms 
\begin{align*}
\widehat{G}(L)[p^n]\simeq \mcal{M}_{G,n}^D(L) 
& \simeq \mrm{Hom}_{G_L}(\mcal{M}_{G,n}(\overline{K}), \mu_{p^n}(\overline{K}))\\
& \simeq \mrm{Hom}_{G_L}(T_p/p^nT_p, \mu_{p^n}(\overline{K}))
\end{align*}
of $G_K$-modules.
Denote by $F$ the action on $T_p$ induced by the $q$-power Frobenius
where $q$ is the order of $k$.
Since $L/K$ is totally ramified and $L\supset K(\mu_{p^{\infty}})$,  
a homomorphism $f\colon T_p/p^nT_p\to \mu_{p^n}(\overline{K})$ commutes with $G_L$
if and only if $f((F-1)x)=0$ for any $x\in T_p/p^nT_p$.
Hence we have natural isomorphisms 
\begin{align*}
\widehat{G}(L)[p^n]
& \simeq \mrm{Hom}_{\mbb{Z}_p}(T_p/((F-1)T_p+p^nT_p), \mu_{p^n}(\overline{K}))\\
& \simeq \mrm{Hom}_{\mbb{Z}_p}(T_p/(F-1)T_p, \mu_{p^n}(\overline{K}))
\end{align*}
of $G_K$-modules.
Taking direct limits over $n$, we obtain an isomorphism of groups 
\begin{equation}
\label{Smith0}
\widehat{G}(L)[p^\infty]
\simeq \mrm{Hom}_{\mrm{cont}}(T_p/(F-1)T_p, \mbb{Q}_p/\mbb{Z}_p).
\end{equation}

Now we denote by $d$ the $\mbb{Z}_p$-rank of $T_p(\mcal{M}_G)$,
and suppose that
the Smith normal form of $F-\mrm{id}$ on $T_p$ is
$$
\mrm{diag}(p^{a_1},\ldots,p^{a_r},0,\ldots,0).
$$
Then 
\begin{equation}\label{Smith1.5}
T_p/(F-1)T_p
\simeq
\bigoplus_{i=1}^r\mbb Z/p^{a_i}\mbb Z
\oplus
\mbb Z_p^{\,d-r}.
\end{equation}
On the other hand,  using the canonical identification 
$p^{-n}T_p/T_p\simeq  \mcal{M}_{G,n}(\bar k)$
of $G_K$-modules,
the same Smith normal form gives
\begin{align}
\label{Smith1}
\mcal{M}_{G,n}(k)
& \simeq \ker (F-\mrm{id} \ \mrm{on}\ p^{-n}T_p/T_p)\notag \\
& \simeq \oplus^r_{i=1} (p^{-a_i}\mbb{Z}_p/\mbb{Z}_p)[p^n]
\oplus (\mbb{Q}_p/\mbb{Z}_p)[p^n]^{d-r}. 
\end{align}
Fix a $\mbb{Z}_p$-basis of $T_p$ and choose
$U,V\in GL_d(\mbb{Z}_p)$ such that
$$
U(F-\mrm{id})V
=\mrm{diag}(p^{a_1},\ldots,p^{a_r},0,\ldots,0).
$$
The automorphisms $U$ and $V$ of $T_p$
induce automorphisms of $p^{-n}T_p/T_p$
for every $n$, compatible with the natural
transition maps in $n$.
Consequently, the isomorphisms in
\eqref{Smith1} can be chosen compatibly
with these transition maps.
Therefore, taking direct limits in~\eqref{Smith1}, 
we obtain a
non-canonical isomorphism
\begin{equation}
\label{Smith2}
\varinjlim_n \ \mcal{M}_{G,n}(k)
\simeq 
\oplus^r_{i=1} (p^{-a_i}\mbb{Z}_p/\mbb{Z}_p)\oplus 
(\mbb{Q}_p/\mbb{Z}_p)^{d-r}.
\end{equation}
On the other hand, it follows from \eqref{Smith0} and \eqref{Smith1.5}
that we have a non-canonical isomorphism of groups
\begin{align}\label{Smith3}
\widehat{G}(L)[p^{\infty}]
& \simeq  \mrm{Hom}_{\mrm{cont}}(\oplus^r_{i=1} (\mbb{Z}_p/p^{a_i}\mbb{Z}_p)\oplus 
\mbb{Z}_p^{d-r},\mbb{Q}_p/\mbb{Z}_p) \notag \\
& \simeq \oplus^r_{i=1} (p^{-a_i}\mbb{Z}_p/\mbb{Z}_p)\oplus (\mbb{Q}_p/\mbb{Z}_p)^{d-r}.
\end{align}
Comparing \eqref{Smith2} and \eqref{Smith3} proves the assertion.
\end{proof}

\begin{theorem}\label{hatG_LD}
Let $K$ be a $p$-adic field and 
$L/K$ a totally ramified Galois extension
contining $K(\mu_{p^{\infty}})$.
Let $A_{/K}$ have good reduction.
Assume that both {\rm (C1)} and {\rm (C2)} hold for $A_{/K}$ and $L/K$.
Then, there exists a non-canonical isomorphism of groups 
$$
A(L)_{\mrm{tor}}
\simeq 
\bar{A}^{\vee}(k)[p^{\infty}]\oplus A(K)_{\mrm{tor}}.
$$
\end{theorem}

\begin{proof}
The assertion follows from
Theorem \ref{LD:p-adic} and Corollary \ref{hatG_p-adic}.
\end{proof}

\subsection{Formal torsion via Smith invariants}
\label{FG}

In this section, we assume that $K$ is a {\it $p$-adic field},
that is, a finite extension of $\mbb{Q}_p$.
We denote by $q=p^f$ the order of its residue field $k$.
The aim of this section is to determine the order and group structure of
the formal part $\widehat{G}(L)[p^\infty]$ in terms of the Smith invariants
of a certain Frobenius action
under (C2) and the assumption that $L/K$ is a totally ramifeid Galois extension.

\subsubsection{Admissible uniformizers}
For any $p$-adic field $K$, we put
$U_K^{(0)}=\mathcal{O}_K^\times$ and
$U_K^{(n)}=1+\mathfrak{m}_K^n$ for $n\geq 1$,
where $\mathfrak{m}_K$ is the maximal ideal of
$\mathcal{O}_K$.
Let $L/K$ be a totally ramified abelian extension.
Let $\mrm{Art}_K:K^{\times} \to \mrm{Gal}(K^{\mrm{ab}}/K)$ be the local Artin map
and denote by 
$\mrm{Art}_{L/K}\colon K^{\times }\to \mrm{Gal}(L/K)$ the composite of 
$\mrm{Art}_K$ and the restriction map 
$\mrm{Gal}(K^{\mrm{ab}}/K)\to \mrm{Gal}(L/K)$.
Here we normalize the local Artin map so that a uniformizer induces arithmetic Frobenius on the maximal unramified extension.
We set
$$
H_L:=\ker \mrm{Art}_{L/K}\subset K^{\times}
$$
and 
$$
b_L:=\inf \left\{ v_p(\mrm{Nr}_{K/\mbb{Q}_p}(u)-1) 
\mid u\in H_L\cap \cO_K^{\times} \right\}
\in \mbb{Z}_{\ge 0}\cup \{\infty\}.
$$

\begin{proposition}
\label{adm}
There exists a uniformizer $\pi$ of $K$
such that $H_L=\pi^{\mbb{Z}}\cdot (H_L\cap \cO_K^{\times})$.
\end{proposition}
\begin{proof}
Take any uniformizer $\pi_0$ of $K$.
Since $L/K$ is totally ramified, 
we have $\mrm{Art}_{L/K}(\cO_K^{\times})=\mrm{Gal}(L/K)$.
Thus there exists a unit $u\in \cO_K^{\times}$ such that 
$\mrm{Art}_{L/K}(\pi_0)^{-1}=\mrm{Art}_{L/K}(u)$.
Putting $\pi=\pi_0u$, we have $\pi\in H_L$. 
The result immediately follows.
\end{proof}

\begin{definition}
Let $L/K$ be a totally ramified abelian extension.
We call a uniformizer $\pi$ of $K$ satisfying $H_L=\pi^{\mbb{Z}}\cdot (H_L\cap \cO_K^{\times})$ 
{\it admissible for $L/K$}.
This is equivalent to the condition that 
$\mrm{Art}_K(\pi)|_L=\mrm{id}_L$.
\end{definition}

\begin{example}
\label{b:LT}
Let $K/J$ be a finite unramified extension
and let $\pi$ be a uniformizer of $J$.
We denote by $J_{\pi}/J$ the Lubin-Tate extension associated with $\pi$,
and, for $n\ge 0$, denote by $J_{\pi,n}$ the $n$-th layer of the Lubin-Tate extension $J_{\pi}$.
A typical example is the case where $J=\mbb{Q}_p$ and $\pi=p$; in this case, $J_{\pi,n}=\mbb{Q}_p(\mu_{p^{n}})$
for $0\le n\le \infty$.
\begin{itemize}
\item[(1)] The field $KJ_{\pi}$ is totally ramified and abelian over $K$.
We have 
$$
H_{KJ_{\pi}}=\mrm{Nr}_{K/J}^{-1}(\pi^{\mbb{Z}})
=\pi^{\mbb{Z}}\cdot \ker \left(\mrm{Nr}_{K/J}\colon \cO_K^{\times}\to \cO_J^{\times} \right).
$$
Thus $\pi$ is admissible for $KJ_{\pi}/K$
and 
$$
b_{KJ_{\pi}}=\infty.
$$ 
Furthermore, 
since $H_{KJ_{\pi}}\subset H_{KJ_{\pi,n}}$,
$\pi$ is admissible also for $KJ_{\pi,n}/K$.
In particular, taking $J=\mbb{Q}_p$ and $\pi=p$, 
we may choose $p$ as an admissible uniformizer
for $K(\mu_{p^{\infty}})/K$
and $K(\mu_{p^n})/K$.

\item[(2)] Suppose $J=K$ and $K/\mbb{Q}_p$ is unramified.
Then $\pi$ is admissible for $K_{\pi}/K$.
We have 
\begin{align*}
b_{K_{\pi,n}}= 
\left\{
\begin{array}{cl}
n
&\quad 
\mbox{$p>2$, or $p=2$ and $n>0$}, \cr
1 
&\quad \mbox{$p=2$ and $n=0$}
\end{array}
\right.
\end{align*}
since 
$H_{K_{\pi,n}}=\pi^\mbb{Z}\cdot U_K^{(n)}$
and $\mrm{Nr}_{K/\mbb{Q}_p}(U_K^{(n)})=U_{\mbb{Q}_p}^{(n)}$.
\end{itemize}
\end{example}

\subsubsection{Smith invariants}

Throughout this subsection, we assume that
$L/K$ is a totally ramified Galois extension
and that (C2) holds for $A_{/K}$ and $L/K$.

Put $\mcal{M}_{G,n}=(\mcal{G}[p^n]^D)^{\aet}$
and denote by $\mcal{M}_G$ the \'etale 
$p$-divisible group associated with $(\mcal{M}_{G,n})_n$.
If $A$ has good reduction, 
then $G=A$, and we have natural
isomorphisms
$\mcal{M}_{G,n}(\overline K)\simeq \bar{A}^{\vee}(\bar k)[p^n]$ 
and $T_p(\mcal{M}_G)\simeq T_p(\bar{A}^{\vee})$.
Let $\mcal{M}^D_{G,n}$ and $\mcal{M}^D_G$
be the Cartier duals of $\mcal{M}_{G,n}$ and $\mcal{M}_G$, respectively.
By Cartier duality, we have   
$T_p(\mcal{M}^D_G)\simeq T_p(\mcal{M}_G)^{\vee}(1)$.

Let $\bar{A}^{\vee}$  denote the dual abelian variety of $\bar{A}$.
This is canonically isomorphic to the reduction of the dual abelian variety 
$A^{\vee}$ of $A$. 

\begin{lemma} 
\label{LtoM}
Let $M$ be the maximal abelian subextension of $L/K$.
Then, there is a natural isomorphism
$$
\widehat{G}(L)[p^{\infty}]
=
\widehat{G}(M)[p^{\infty}]
\simeq (T_p(\mcal{M}^D_G)\otimes_{\mbb{Z}_p}\mbb{Q}_p/\mbb{Z}_p)^{G_M}
$$
of $G_K$-modules.
\end{lemma}

\begin{proof}
By Corollary \ref{hatG_cyclic},
$\mrm{Gal}(L/M)$ acts trivially on $\widehat{G}(L)[p^{\infty}]$,
that is,
$\widehat{G}(L)[p^{\infty}]=\widehat{G}(M)[p^{\infty}]$.
By Theorem \ref{hatG},
 for every $n>0$ there are natural
isomorphisms
$$
\hat{G}(L)[p^n]\simeq \mcal{M}_{G,n}^D(L)
\simeq 
(T_p(\mcal{M}_G^D)\otimes_{\mbb{Z}_p} p^{-n}\mbb{Z}_p/\mbb{Z}_p)^{G_L}. 
$$
Taking $G_M$-fixed parts, 
$$
\hat{G}(M)[p^n]\simeq 
(T_p(\mcal{M}_G^D)\otimes_{\mbb{Z}_p} p^{-n}\mbb{Z}_p/\mbb{Z}_p)^{G_M}. 
$$
These identifications are compatible with the transition maps in $n$.
Taking direct limits over $n$ gives the desired isomorphism.
We remark that, since every element on the right has finite $p$-power order,
taking $G_M$-invariants commutes with this direct limit.
\end{proof}

We denote by 
$$
\rho_G \colon G_K\to GL_{\mbb{Z}_p}(T_p(\mcal{M}_G))
$$ 
the representation induced by the natural action
of $G_K$ on $T_p(\mcal{M}_G)$.
Since $\mcal{M}_G$ is \'etale, $\rho_G$ is unramified.
In particular, $\rho_G$ has abelian image.
By abuse of notation,
we use the same notation $\rho_G$ for the 
induced representation of $\mrm{Gal}(K^{\mrm{ab}}/K)$;
$$
\rho_G\colon \mrm{Gal}(K^{\mrm{ab}}/K)\to GL_{\mbb{Z}_p}(T_p(\mcal{M}_G)).
$$
Let $M$ be the maximal abelian subextension of $L/K$, 
and let $\pi$ be an admissible uniformizer of $M/K$.
Set 
$$
\Fr_G:=\rho_G(\mrm{Art}_K(\pi)).
$$
By our normalization of the Artin map, 
this coincides with the image of the $q$-th arithmetic Frobenius of 
the unramified representation $\rho_G$.
We also set 
$$
D_{K,\pi}:=\ker \left( \frac{q}{\mrm{Nr}_{K/\mbb{Q}_p}(\pi)} - \Fr_G \ \mrm{on}\ T_p(\mcal{M}_G)\otimes_{\mbb{Z}_p} \mbb{Q}_p/\mbb{Z}_p\right).
$$
For example, if $A$ has good reduction and $q=\mrm{Nr}_{K/\mbb{Q}_p}(\pi)$,
then 
$$
D_{K,\pi}\simeq \ker \left( 1 - \varPhi \ \mrm{on}\ \bar{A}^{\vee}(\bar k)[p^{\infty}]\right)
\simeq \bar{A}^{\vee}(k)[p^{\infty}]
$$
where $\varPhi$ is the Frobenius endomorphism of $\bar{A}^{\vee}$.
\medskip

Let us state our main result of this section.
To simplify notation, we put 
\begin{align*}
C_{p^n}:=
\left\{
\begin{array}{cl}
\mbb{Z}/p^n\mbb{Z} 
&\quad 
\mrm{if}\ n\in \mbb{Z}_{\ge 0}, \cr
\mbb{Q}_p/\mbb{Z}_p 
&\quad \mrm{if}\ n=\infty.
\end{array}
\right.
\end{align*}

\begin{theorem}
\label{FG:Main}
Assume that (C2) holds for $A_{/K}$ and $L/K$.
Assume that  $L/K$ is totally ramified and
let $M$ be the maximal abelian subextension of $L/K$.
Let $\pi$ be an admissible uniformizer for $M/K$.
Write $g=\dim A$, $t=\mrm{rank}\,T$ for the toric rank of $A$ and
$f_B$ for the $p$-rank of $\bar{B}$.
Let  $b_M$ and $D_{K,\pi}$ be as above.

Suppose that  the Smith normal form of 
$q\mrm{Nr}_{K/\mbb{Q}_p}(\pi)^{-1} - \Fr_G$
on the $\mbb{Z}_p$-module $T_p(\mcal{M}_G)$ is
$$
\operatorname{diag}(p^{a_1},\dots ,p^{a_r},0,\dots ,0)
$$
where $0\le a_1\le \cdots \le a_r$ and $0\le r\le t+f_B$.
\begin{itemize}
\item[{\rm (1)}] We have 
\begin{equation}
\label{formalstructure}
\widehat{G}(L)[p^{\infty}]\simeq D_{K,\pi}[p^{b_M}]\simeq 
\bigoplus^r_{i=1} C_{p^{\mrm{min}\{a_i,b_M\}}}\oplus (C_{p^{b_M}})^{t+f_B-r}
\end{equation}
\item[{\rm (2)}] Assume moreover that $A$ has good reduction,
and let $f_A$ denote the $p$-rank of
$\bar A$. 
Then $\det (q\mrm{Nr}_{K/\mbb{Q}_p}(\pi)^{-1} - \Fr_A)\not=0$ if and only if 
$P_{A,K}(\mrm{Nr}_{K/\mbb{Q}_p}(\pi))\not=0$.
If these equivalent conditions hold, then $r=f_A$ and 
\begin{equation}
\label{formalorder}
v_p(\det (q\mrm{Nr}_{K/\mbb{Q}_p}(\pi)^{-1} - \Fr_A))
=\sum^{f_A}_{i=1} a_i=v_p(P_{A,K}(\mrm{Nr}_{K/\mbb{Q}_p}(\pi)))-fg.
\end{equation}
In particular,  
$$
\#\widehat{A}(L)[p^{\infty}]\le p^{v_p(P_{A,K}(\mrm{Nr}_{K/\mbb{Q}_p}(\pi)))-fg}.
$$
Moreover, if $f_A>0$, equality holds if and only if
$b_M\ge a_{f_A}$.
If $f_A=0$, equality always holds.
\end{itemize}
\end{theorem}

\begin{proof}
Let $\rho^D_G$ denote the representation of $\mrm{Gal}(K^{\mrm{ab}}/K)$
acting on $T_p(\mcal{M}^D_G)$;
$$
\rho^D_G\colon \mrm{Gal}(K^{\mrm{ab}}/K)\to GL_{\mbb{Z}_p}(T_p(\mcal{M}^D_G)).
$$
Let $U_{K,\pi}:=\rho^D_G(\mrm{Art}_K(\pi))$, and put 
$$
D_{K,\pi}':=\ker \left(\rho^D_G(\mrm{Art}_K(\pi)) - \mrm{id} \ \mrm{on}\ T_p(\mcal{M}^D_G)\otimes_{\mbb{Z}_p} \mbb{Q}_p/\mbb{Z}_p\right)
$$
Choose a $\mbb Z_p$-basis of $T_p(\mcal M_G)$ and the dual basis of
$T_p(\mcal M_G^D)$ under the Cartier duality pairing.
Regarding $\Fr_G$ and $U_{K,\pi}$ as the matrices of the corresponding
actions with respect to these bases, Cartier duality gives
$
U_{K,\pi}=q\cdot \mrm{Nr}_{K/\mbb{Q}_p}(\pi)^{-1}({}^t\Fr_G)^{-1}.
$
Hence 
$$
U_{K,\pi}-I
=
({}^t\Fr_G)^{-1} 
\cdot {}^t\left(
q\,\mrm{Nr}_{K/\mbb Q_p}(\pi)^{-1}I-\Fr_G
\right).
$$
Since multiplication by an invertible matrix and transposition do not
change the Smith invariants, $U_{K,\pi}-I$ and
$q\mrm{Nr}_{K/\mbb Q_p}(\pi)^{-1}I-\Fr_G$
have the same Smith invariants.
Thus there is a non-canonical isomorphism of groups
$$
D_{K,\pi}\simeq D'_{K,\pi}.
$$

\medskip

\noindent
(1) Let
$P\in T_p(\mcal M_G^D)\otimes_{\mbb Z_p}\mbb Q_p/\mbb Z_p.
$
By Lemma~\ref{LtoM}, $P$ belongs to
$\widehat G(L)[p^\infty]=\widehat G(M)[p^\infty]$
if and only if
$$
\rho_G^D(\mrm{Art}_K(x))(P)=P
\qquad (x\in H_M).
$$
By Proposition \ref{adm}, 
this is equivalent to the following two conditions. 
\begin{itemize}
\item[(i)]  $\rho^D_G(\mrm{Art}_K(\pi))(P)=P$, and 
\item[(ii)] $\rho^D_G(\mrm{Art}_K(u))(P)=P$ for every $u\in H_M\cap \cO_K^{\times}$.
\end{itemize}
Condition (i) is equivalent to 
$P\in D_{K,\pi}'$.
Since $\rho_G^D$ restricted to $I_K$ is given by $\chi_p$, 
local class field theory gives  
$\rho^D_G\circ \mrm{Art}_K(u)=\mrm{Nr}_{K/\mbb{Q}_p}(u)^{-1}$
as a scalar for $u\in \cO_K^{\times}$.
Thus, by definition of $b_M$, condition (ii) 
is equivalent to $p^{b_M}P=0$.
Therefore, we have
$\widehat{G}(L)[p^{\infty}]
=D_{K,\pi}'\cap \widehat{G}(\overline{K})[p^{b_M}]
=D_{K,\pi}'[p^{b_M}]\simeq D_{K,\pi}[p^{b_M}]$.
The assertion now follows from the elementary divisor theory 
applied to $q\mrm{Nr}_{K/\mbb{Q}_p}(\pi)^{-1} - \Fr_G$.

\medskip

\noindent
(2) First note that $t=0$, $f_B=f_A$ and 
$T_p(\mcal{M_G})=T_p(\bar{A}^{\vee})$, 
since $A$ has good reduction.
Let $\beta_1,\dots \beta_{f_A}$ be the eigenvalues of $\Fr_A$.
By a result of Manin, the characteristic polynomial
$P_{A,K}(T)$ of the $q$-power Frobenius acting on $\bar A$
is also the characteristic polynomial of the Frobenius acting on the
Dieudonn\'e module of $\bar A^\vee[p^\infty]$;
see \cite[p.~96, Corollary]{Dem72}.
Moreover, its $p$-adic unit roots correspond precisely to the slope-zero
part of the Dieudonn\'e module
\cite[p.~98, Theorem]{Dem72},
which is identified with the \'etale Tate module
$T_p(\bar A^\vee)$ by \cite{Ill79}.
Thus $\beta_1,\dots,\beta_{f_A}$ are precisely the unit roots of
$P_{A,K}(T)$.
See also the proof of \cite[Lemma~3.9]{Cre26}.
Write $\alpha_1,\ldots,\alpha_{2g-f_A}$ for the remaining roots of
$P_{A,K}(T)$.
Then
$$
P_{A,K}(q\mrm{Nr}_{K/\mbb{Q}_p}(\pi)^{-1})=
\prod^{f_A}_{i=1}\left(q\mrm{Nr}_{K/\mbb{Q}_p}(\pi)^{-1}-\beta_i\right)
\cdot 
\prod^{2g-f_A}_{i=1}\left(q\mrm{Nr}_{K/\mbb{Q}_p}(\pi)^{-1}-\alpha_i\right).
$$
Since each $q\mrm{Nr}_{K/\mbb{Q}_p}(\pi)^{-1}-\alpha_i$ is a $p$-adic unit,
we find 
\begin{align}\label{Sm1}
v_p\left(P_{A,K}(q\mrm{Nr}_{K/\mbb{Q}_p}(\pi)^{-1})\right)
& =v_p\left(\prod^{f_A}_{i=1}
(q\mrm{Nr}_{K/\mbb{Q}_p}(\pi)^{-1}-\beta_i)\right) \notag\\
&=v_p\left(\det\left( q\mrm{Nr}_{K/\mbb{Q}_p}(\pi)^{-1}-\Fr_A \right) \right).
\end{align}
On the other hand, by the functional equation
$q^{g}P_{A,K}(T)=T^{2g}P_{A,K}(q/T)$ of the Weil polynomial,
we have 
\begin{align}\label{Sm2}
v_p\left(P_{A,K}(q\mrm{Nr}_{K/\mbb{Q}_p}(\pi)^{-1})\right)
&=v_p\left(P_{A,K}(\mrm{Nr}_{K/\mbb{Q}_p}(\pi))\cdot 
q^g\mrm{Nr}_{K/\mbb{Q}_p}(\pi)^{-2g}\right) \notag \\
&=v_p\left(P_{A,K}(\mrm{Nr}_{K/\mbb{Q}_p}(\pi))\right)-fg.
\end{align}
Therefore, it follows from \eqref{Sm1} and \eqref{Sm2} that 
$$
v_p\left(\det\left( q\mrm{Nr}_{K/\mbb{Q}_p}(\pi)^{-1}-\Fr_A \right) \right)
=v_p\left(P_{A,K}(\mrm{Nr}_{K/\mbb{Q}_p}(\pi))\right)-fg.
$$
The remaining assertions follow immediately from {\rm (1)} and the
Smith normal form.
\end{proof}

\section{Applications and higher-dimensional examples}\label{sec:applications}

\subsection{Torsion over pro-$p$ and cyclotomic extensions}\label{subsec:cyclotomic}

We begin with a general rigidity statement for torsion over pro-$p$
extensions. Throughout this section, as before, 
$p$ is an odd prime unless otherwise specified.

\begin{proposition}\label{prop:global-prop-rigidity}
Let $F_\infty$ be a pro-$p$ Galois extension of a number field $F$,
let $v\mid p$ be a finite place of $F$, and fix a place $w$ of
$F_\infty$ above $v$.
Let $D_w\subset\Gal(F_\infty/F)$ be the decomposition subgroup at $w$,
and let $M/F_v$ be the corresponding local extension.
Let $A_{/F}$ be an abelian variety.
Assume that $M/F_v$ is totally ramified and that
$A_{/F_v}$ acquires semistable reduction over a finite extension
$K'/F_v$ satisfying $e_{K'}<p-1$.
Then
\[
 A(F_\infty)_{\rm tor}
 =
 A(F_\infty^{D_w})_{\rm tor}.
\]
In particular, if $w$ is the unique place of $F_\infty$ above $v$, then
\[
 A(F_\infty)_{\rm tor}=A(F)_{\rm tor}.
\]
\end{proposition}

\begin{proof}
Take any $P\in A(F_\infty)_{\rm tor}$.
Via the embedding of $F_\infty$ into $M$ determined by $w$,
we may regard $P$ as an element of $A(M)_{\rm tor}$.
Since $M/F_v$ is totally ramified and pro-$p$,
Theorem~\ref{Yoshida_descent} gives
$P\in A(F_v)_{\rm tor}$.
Hence $P$ is fixed by $D_w$, and therefore
$P\in A(F_\infty^{D_w})_{\rm tor}$.
This proves the first assertion.
If $w$ is the unique place above $v$, then
$D_w=\Gal(F_\infty/F)$, so that
$F_\infty^{D_w}=F$.
The second assertion follows.
\end{proof}

We next apply this rigidity statement to cyclotomic $\mbb Z_p$-extensions.
For this purpose, we use a uniform bound for the ramification degree
required to acquire semistable reduction.
For any $g\ge 1$, put 
$$
m_g=\mrm{max}\left\{\mrm{lcm}(n_1,\dots ,n_r)
\mid r>0,\ n_1,\dots, n_r>0,\  \sum^r_{i=1}\varphi(n_i)\le 2g\right\}.
$$
\begin{table}[htbp]
\centering
\setlength{\tabcolsep}{4pt}
\renewcommand{\arraystretch}{1.2}
\begin{tabular}{c|rrrrrrrrrr}
\hline
$g$
& 1 & 2 & 3 & 4 & 5 & 6 & 7 & 8 & 9 & 10 \\
\hline
$m_g$
& 6 & 12 & 30 & 60 & 120 & 210 & 420
& 840 & 1260 & 2520 \\
$p_g$
& 11 & 17 & 37 & 67 & 127 & 223 & 431
& 853 & 1277 & 2531 \\
\hline
\end{tabular}
\caption{Values of $m_g$ and $p_g$,
where $p_g$ is the smallest prime
strictly greater than $m_g+1$.}
\label{tab:mg}
\end{table}

\begin{corollary}\label{cor:cyclotomic-rigidity}
Let $F$ be a number field and put
$
F_{\infty,p}:=F\mathbb Q_{\infty,p}.
$
Assume that there exists a finite  place $v$ of $F$ above $p$
which is unramified over $\mathbb Q$.
Let $A_{/F}$ be a $g$-dimensional abelian variety.
\begin{itemize}
\item[{\rm (1)}]
If $A$ has semistable reduction at $v$, then
$$
A(F_{\infty,p})_{\rm tor}=A(F)_{\rm tor}.
$$

\item[{\rm (2)}]
If $p>m_g+1$, then
$$
A(F_{\infty,p})_{\rm tor}=A(F)_{\rm tor}.
$$
\end{itemize}
\end{corollary}

\begin{lemma}\label{semistabledegree}
Assume $p>2g+1$.
Let $A$ be a $g$-dimensional abelian variety over a $p$-adic field $K$.
Then, $A$ acquires semistable reduction over a finite extension $K'$ of $K$ 
such that  $e_{K'/K}\le m_g$.
\end{lemma}

\begin{proof}
For any prime $\ell\not=p$,
let $I_{A,K}$ be the subgroup of $I_K$
consisting of those $\sigma\in I_K$ which act unipotently on
$T_\ell(A)$.
It is known that $I_{A,K}$ is independent of $\ell$ 
and is an open subgroup of $I_K$.
Moreover, $A$ has semistable reduction over $K'$ if and only if $I_{K'}\subset I_{A,K}$
(cf.\ \cite[Theorem 4.2]{SiZa98}).
We denote by 
$
\varPhi_{A,K}=I_K/I_{A,K}
$
the local monodromy group of $A_{/K}$,
and set $m:=\# \varPhi_{A,K}$.
It suffices to show 
$$
m\le m_g.
$$
Indeed, once this is proved, since $\varPhi_{A,K}$ is a quotient of the
tame inertia group, one can choose a finite tamely ramified extension
$K'/K$ with
$$
e_{K'/K}\le m
$$
such that $I_{K'}\subset I_{A,K}$.

By $p>2g+1$ and \cite[Corollary 6.3]{SiZa98}, we have $p\nmid m$.
Hence the wild inertia subgroup $P_K$ acts trivially on
$\varPhi_{A,K}$, and therefore there is a surjection
$I_K/P_K\to \varPhi_{A,K}$.
It follows that $\varPhi_{A,K}$ is cyclic of order $m$.
Let $t_v$ and $a_v$ denote the toric rank and the abelian rank,
respectively, of the special fiber of the N\'eron model of $A/K$.
Take a finite extension over which $A$ has semistable reduction,
and let $t$ and $a$ denote the corresponding toric and abelian ranks.
By \cite[Theorem 5.2]{SiZa98},
there exists an injective representation
$$
\rho\colon \varPhi_{A,K}\to GL_{t-t_v}(\mbb{Z})\times Sp_{2(a-a_v)}(\mbb{Q}_{\ell})
$$
such that the projection onto the first factor is independent of $\ell\not=p$,
and the characteristic polynomial of the projection onto the second factor has integer coefficients independent of $\ell\not=p$.
Thus we may regard $\rho$ as an injective representation
$
\rho\colon
\varPhi_{A,K}
\hookrightarrow
GL_{t-t_v+2(a-a_v)}(\mbb Q_\ell),
$
whose characteristic polynomials have rational coefficients.
Fix a generator $\gamma$ of $\varPhi_{A,K}$.
Since $\rho(\gamma)$ is of finite order $m$,
$\rho(\gamma)$ is diagonalizable over $\overline{\mbb{Q}}_{\ell}$
and all of its eigenvalues are roots of unity.
Since the characteristic polynomial $P(T)$ of 
$\rho(\gamma)$ has rational coefficients,
 we may write
$$
P(T)=\prod_{i=1}^r\varPhi_{n_i}(T)^{a_i},
\qquad
m=\mrm{lcm}(n_1,\ldots,n_r),
$$
for some integers $n_i>0$ and $a_i>0$, where
$\varPhi_{n_i}(T)$ denotes the $n_i$-th cyclotomic polynomial.
Since $\deg P(T)=t-t_v+2(a-a_v)\le 2g$,
we have 
$$
\sum^r_{i=1} \phi(n_i)\le \sum^r_{i=1}a_i\vphi(n_i)=\deg P(T)\le 2g. 
$$
Therefore, we have $m\le m_g$.
This completes the proof.
\end{proof}

\begin{proof}[{Proof of Corollary \ref{cor:cyclotomic-rigidity}}]
Let $v$ be a finite place of $F$ above $p$ such that
$F_v/\mbb Q_p$ is unramified.
The extension $F_{\infty,p}/F$ is pro-$p$, and there is a unique place
$w$ of $F_{\infty,p}$ above $v$.
Moreover, the corresponding local extension over $F_v$ is totally ramified.
Assertion (1) follows immediately from Proposition \ref{prop:global-prop-rigidity}.
We prove (2). Suppose $p>m_g+1$.
Then, we have $p>2g+1$ since $m_g>2g$.
By Lemma~\ref{semistabledegree}, the abelian variety $A_{/F_v}$
acquires semistable reduction over a finite extension $K'/F_v$ such that
$e(K'/F_v)\le m_g$.
Since $F_v/\mbb Q_p$ is unramified,
$$
e_{K'}=e(K'/F_v)\le m_g<p-1.
$$
Hence (2) also follows from
Proposition~\ref{prop:global-prop-rigidity}.
\end{proof}

\begin{remark}
The conclusion of Corollary
\ref{cor:cyclotomic-rigidity} (1)
does not necessarily hold when $p=2$.
Indeed, consider 
the elliptic curve 
$E: y^2+xy+y=x^3-11x+12$ 
(LMFDB label 14.a4).
Since $\Delta_E=98$ and $c_4(E)=505$,
the curve $E$ has multiplicative reduction at $2$.
We also have 
$$
E(\overline{\mbb{Q}})[2]
=\left\{ O, \left(\frac{7}{4},-\frac{11}{8} \right), 
\left(-1+2\sqrt{2}, -\sqrt{2} \right), 
\left(-1-2\sqrt{2}, \sqrt{2} \right) \right\}
$$
Since $\mbb{Q}(\sqrt{2})\subset \mbb{Q}_{\infty,2}$,
we have 
$E(\mbb{Q}_{\infty,2})[2]\not=E(\mbb{Q})[2]$.
\end{remark}

We now apply Proposition \ref{prop:global-prop-rigidity}
to determine the torsion subgroups of modular Jacobians
over certain infinite non-abelian extensions.
We first give a family of such extensions.

\begin{proposition}\label{Free}
Let $p\geq 5$ be a prime, and let $F$ be an imaginary quadratic field.
Assume that $p$ is inert in $F$ and that the class number of $F$
is prime to $p$.
Let $F^{\{p\}}$ denote the maximal pro-$p$ Galois extension of $F$
unramified outside the primes above $p$.
\begin{itemize}
\item[{\rm (1)}]
$\mrm{Gal}(F^{\{p\}}/F)$ is a free pro-$p$ group on two generators.
\item[{\rm (2)}]
Let $v$ be the unique finite place of $F$ above $p$.
Then $F^{\{p\}}/F$ is totally ramified at $v$.
In particular, there exists a unique place of $F^{\{p\}}$ above $v$.
\end{itemize}
\end{proposition}

\begin{proof}
Assertion (1) follows from \cite[Theorem 2.13 (i)]{MaRu08}.
We prove (2). Let $E/F$ be a finite Galois subextension of $F^{\{p\}}/F$
and set $G:=\mrm{Gal}(E/F)$. Take any finite place $w$ of $E$ above $v$ 
and denote by $I$ the inertia subgroup of $G$.
It suffices to prove $I=G$.
Assume $I\not =G$. 
Since $G$ is a finite $p$-group and $I\neq G$,
there exists a maximal subgroup $H$ of $G$ containing $I$.
(Note that every maximal subgroup of a finite $p$-group is normal and has index $p$.)
By construction, the extension $E^H/F$ is cyclic of order $p$.
Furthermore, $E^H/F$ is unramified at $v$, and 
it is also unramified outside $p$ by $F\subset E^H\subset F^{\{p\}}$.
Since $p$ is inert in $F$, $E^H/F$ is unramified at every finite place.
Therefore there exists a nontrivial unramified cyclic extension of $F$
of degree $p$, so $p$ divides the class number of $F$,
contrary to our assumption.
Hence $I=G$ as desired.
\end{proof}

\begin{example}\label{Jac_Ex}
Let $N\ge 5$ be a prime and put $J:=J_0(N)$.
Then $J$ has purely toric reduction at $N$
and good reduction at every prime $q\neq N$;
see \cite[\S2]{Rib99} and the references therein.
By \cite[Theorem 1]{Maz77}, we have 
$J(\mbb{Q})_{\mrm{tor}}\simeq\mbb{Z}/n\mbb{Z}$
where $n$ is the numerator of $(N-1)/12$.

Let $F$ and $v\mid p$ be as in Proposition \ref{Free}.
Since $J$ has semistable reduction at every prime,
Propositions \ref{prop:global-prop-rigidity}
and \ref{Free} imply
$$
J(F^{\{p\}})_{\mrm{tor}}=J(F)_{\mrm{tor}}.
$$
Thus it suffices to determine $J(F)_{\mrm{tor}}$.
Let $q_1,q_2\ge 3$ be distinct primes,
different from $N$ and unramified in $F$,
and let $w_i\mid q_i$ be primes of $F$.
Since $e_{F_{w_i}}=1<q_i-1$, reduction induces injections
$
J(F)_{\mrm{tor}}\hookrightarrow
J(k_{w_i})
$
for each $i$.
Hence, we have 
$$
J(F^{\{p\}})_{\mrm{tor}}
\simeq\mbb{Z}/n\mbb{Z}\quad \mrm{if}\quad 
\gcd(\#J(k_{w_1}),\#J(k_{w_2}))=n.
$$
For instance, take $F=\mbb{Q}(i)$ and $p=7$.
The following table determines
$J(F^{\{7\}})_{\mrm{tor}}$
for every prime level $N$ such that $X_0(N)$ has genus $2$.
\begin{center}
\renewcommand{\arraystretch}{1.23}
\begin{tabular}{ccccc}
\toprule
$N$ &
$\#J(\mbb{F}_9)$ &
$\#J(\mbb{F}_5)$ &
g.c.d. &
$J(F^{\{7\}})_{\mrm{tor}}$\\
\midrule
23 & 121 & 44 & 11 & $\Z/11\Z$\\
29 & 161 & 49 & 7 & $\Z/7\Z$\\
31 & 80 & 25 & 5 & $\Z/5\Z$\\
37 & 105 & 48 & 3 & $\Z/3\Z$\\
\bottomrule
\end{tabular}
\end{center}
The point counts in the table can be verified
by the following MAGMA code, using the
characteristic polynomials of Hecke operators.
{\small
\begin{verbatim}
for N in [23, 29, 31, 37] do
    S := CuspidalSubspace(ModularSymbols(N, 2, 1));
    assert Dimension(S) eq 2;

    h3 := CharacteristicPolynomial(HeckeOperator(S, 3));
    h5 := CharacteristicPolynomial(HeckeOperator(S, 5));

    J9 := Integers()!(Evaluate(h3, 4)*Evaluate(h3, -4));
    J5 := Integers()!Evaluate(h5, 6);

    d := GCD(J9, J5);
    n := Numerator((N-1)/12);
    assert d eq n;

    printf "N=%o, J(F9)=%o, J(F5)=%o, gcd=%o\n",
           N, J9, J5, d;
end for;
\end{verbatim}
}

\end{example}

\medskip

We now turn to the full cyclotomic extension
$F(\mu_{p^\infty})/F$.
The next result gives a strong restriction on the field of definition
of any torsion point which is not already defined over $F$.
It follows from local decomposition and the cyclotomic action on the 
formal torsion.

\begin{theorem}\label{thm:E1}
Let $F$ be a number field, let $p$ be an odd prime, and let
$\pp\mid p$ be a prime of $F$ such that
$F_{\pp}/\mathbb Q_p$ is unramified.
Let $A/F$ be an abelian variety with good reduction at $\pp$.
For every
$P\in A(F(\mu_{p^\infty}))_{\tor}\setminus A(F)_{\tor}$,
there is an integer $r\ge1$ such that
\[
 F(P)=F(\mu_{p^r}).
\]
\end{theorem}

The key input is the following no-growth statement for fixed fields
of nontrivial finite subgroups of the cyclotomic Galois group.

\begin{proposition}\label{prop:cyclotomic-fixed-field}
Let $F$, $p$, $\pp$, and $A$ be as in Theorem~\ref{thm:E1}, and put
$M=F(\mu_{p^\infty})$.
For every nontrivial finite subgroup $\Delta\subset\Gal(M/F)$, we have
\[
 A(M^\Delta)_{\tor}=A(F)_{\tor}.
\]
\end{proposition}

\begin{proof}
Put $K=F_{\pp}$. 
Write $E_n=F(\mu_{p^n})$, $K_n=K(\mu_{p^n})$, 
and $K_\infty=K(\mu_{p^\infty})$.
Since $K/\mathbb Q_p$ is unramified, $K_n/K$ is totally ramified and
\[
 (p-1)p^{n-1}=[K_n:K]\le [E_n:F]\le(p-1)p^{n-1}.
\]
Thus $\pp$ has a unique, totally ramified extension $w$ to $M$,
and localization gives
\[
 D_w(M/F)=I_w(M/F)=\Gal(M/F)
 \simeq\Gal(K_\infty/K)\simeq\mathbb Z_p^\times.
\]
We identify $\Delta$ with its image in $\Gal(K_\infty/K)$.

By Propositions~\ref{pro-nilpotent} and~\ref{C2}(a),
conditions (C1) and (C2) hold for $A/K$ and $K_\infty/K$.
Theorem~\ref{LD:p-adic} therefore gives
\[
 A(K_\infty)_{\tor}
 =
 \widehat A(K_\infty)[p^\infty]\oplus A(K)_{\tor}.
\]
By Corollary~\ref{hatG_cyclic}(2), the action on the first summand
is given by $\chi_p$.
Choose $1\ne\delta\in\Delta$.
Since $\delta$ has finite order,
$\chi_p(\delta)-1\in\mathbb Z_p^\times$.
Hence $\widehat A(K_\infty)[p^\infty]^\Delta=0$, 
and taking $\Delta$-invariants yields
\[
 A(K_\infty^\Delta)_{\tor}=A(K)_{\tor}.
\]
Now let $P\in A(M^\Delta)_{\tor}$.
Its image under localization belongs to $A(K)_{\tor}$.
Since the decomposition group at $w$ is all of $\Gal(M/F)$,
every element of $\Gal(M/F)$ fixes $P$.
Thus $P\in A(F)_{\tor}$, proving the assertion.
\end{proof}

\begin{proof}[Proof of Theorem~\ref{thm:E1}]
Put $M=F(\mu_{p^\infty})$.
Via the cyclotomic character, write
\[
 G:=\Gal(M/F)=T\times\Gamma,\qquad
 \#T=p-1,\quad \Gamma=1+p\mathbb Z_p,
\]
where $T$ is the torsion subgroup of $\mathbb Z_p^\times$.
For $P$ as in the theorem, let $S=\Gal(M/F(P))$, an open subgroup of $G$.
If $S\cap T\ne1$, applying
Proposition~\ref{prop:cyclotomic-fixed-field} with
$\Delta=S\cap T$ would give $P\in A(F)_{\tor}$.
Hence $S\cap T=1$.
The projection $S\to\Gamma$ is injective, so $S$ is pro-$p$ and
its projection to $T$ is trivial.
Consequently $S$ is an open subgroup of $\Gamma$, hence
$S=1+p^r\mathbb Z_p$ for some $r\ge1$.
Its fixed field is $F(\mu_{p^r})$, as required.
\end{proof}

\subsection{Explicit torsion of modular Jacobians over $\mbb{Q}^{\mrm{ab}}$}
\label{Sec_Jac}
Let $N$ be a prime and denote by $n$ the numerator of $(N-1)/12$.
Put $J:=J_0(N)$, the Jacobian of the modular curve $X_0(N)$.
Let $\Sigma$ be the Shimura subgroup of $J$,
namely, the kernel of the natural pullback map
$J_0(N)\to J_1(N)$.
It follows from \cite{Maz77} that
$$
J(\mbb Q)_{\mrm{tor}}\simeq \mbb Z/n\mbb Z
\quad\mrm{and}\quad
\Sigma\simeq\mu_n.
$$
The Jacobian $J$ has purely multiplicative reduction at $N$
and good reduction at every prime $q\neq N$.
For each prime $q\neq N$, we denote by $P_q(T)$ the characteristic
polynomial of the $q$-power Frobenius endomorphism of the reduction
of $J$ at $q$.

The aim of this subsection is to determine the group structure of
$J(\mbb Q^{\mrm{ab}})_{\mrm{tor}}$
for several values of $N$.

\begin{lemma}\label{Jac_p'} 
Assume that there exist distinct primes $q_1,q_2\neq N$
satisfying the following conditions:
\begin{itemize}
\item[{\rm (i)}]
Every prime divisor of $P_{q_1}(1)$ either divides $n$ or is equal to $q_1$.
\item[{\rm (ii)}]
$q_1\nmid P_{q_2}(1)$.
\end{itemize}
Then the prime-to-$n$ part 
of $J(\mbb Q^{\mrm{ab}})_{\mrm{tor}}$ is trivial.
\end{lemma}

\begin{lemma}\label{Jac_p} 
Let $\ell$ be an odd prime divisor of $n$. 
Assume that there exists a prime $q_3\neq \ell,N$
satisfying 
$(X-1)^2\nmid P_{q_3}(X)$
in $\mbb{F}_{\ell}[X]$.
Then, 
$$
J(\mbb Q^{\mrm{ab}})[\ell]= J(\mbb{Q})[\ell]
\oplus \Sigma(\overline{\mbb{Q}})[\ell].
$$
\end{lemma}

The two lemmas are proved simultaneously.

\begin{proof}[Proof of Lemmas \ref{Jac_p'} and \ref{Jac_p}]
For any prime $\ell$,  
put $W_{\ell}=J(\mbb{Q}^{\mrm{ab}})[\ell]$.
Assume that $W_{\ell}\not=0$.
Set $d=\dim_{\mbb{F}_{\ell}} W_{\ell}$.
Since the action of $G_{\mbb Q}$ on $W_\ell$ factors through an abelian
quotient, we have
$(W_{\ell}\otimes_{\mbb{F}_{\ell}} \bar{\mbb{F}}_{\ell})^{\mrm{ss}}
\simeq \oplus^d_{i=1} \bar{\mbb{F}}_{\ell}(\psi_i)$
where $\psi_i\colon G_{\mbb{Q}}\to \bar{\mbb{F}}_{\ell}^{\times}$
is a continuous character.
Since $J$ has semistable reduction everywhere,
each $\psi_i\colon G_{\mbb{Q}}\to \bar{\mbb{F}}_{\ell}^{\times}$
is unramified outside $\ell$.
Indeed, at every prime different from $\ell$ the inertia action on
$J[\ell]$ is unipotent, and hence acts trivially on every
one-dimensional Jordan--Hölder constituent.
By the Kronecker--Weber theorem, $\psi_i$ therefore factors through
$\mrm{Gal}(\mbb{Q}(\mu_{\ell^{\infty}})/\mbb{Q})$.
Since the image of $\psi_i$ has order prime to $\ell$, it follows that
$\psi_i$ factors through
$\mrm{Gal}(\mbb Q(\mu_\ell)/\mbb Q)$ when $\ell>2$.
This shows that 
$\psi_i=\bar{\chi}_{\ell}^{a_i}$ for some $0\le a_i<\ell-1$.
Here $\bar{\chi}_\ell$ denotes the
mod-$\ell$ cyclotomic character.
For $\ell=2$, the same conclusion holds with $a_i=0$, since every finite
abelian extension of $\mbb Q$ unramified outside $2$ has $2$-power
degree, whereas the image of $\psi_i$ has odd order.
For $\ell>2$,  
Caruso's theorem \cite[Th\'eor\`em 1.0.3]{Car06} on tame inertia weights implies that
$a_i\in \{0,1\}$ 
(here we recall that
$J$ has good reduction outside $N$ and has semistable reduction at $N$). 
Consequently, for every prime $\ell$, we have
\begin{equation}\label{Caruso}
W_\ell^{\mrm{ss}}
\simeq
\bigoplus_{i=1}^d
\mbb F_\ell(\overline\chi_{\ell}^{a_i}),
\qquad
a_i\in \{0,1\}.
\end{equation}

Let $q$ be a prime different from $\ell$ and $N$.
Since $W_\ell$ is a $G_{\mbb Q}$-stable subspace of $J[\ell]$, the
characteristic polynomial
$
\det(X-\mrm{Fr}_q\mid W_{\ell})
$
divides $P_q(X)$ modulo $\ell$,
where $\mrm{Fr}_q$ is the arithmetic $q$-Frobenius. 
By \eqref{Caruso}, every root of this characteristic polynomial is
either $1$ or $q$.
Since $W_{\ell}\neq0$, we therefore have
$P_q(1)\equiv 0$ mod $\ell$ 
or $P_q(q)\equiv 0$ mod $\ell$.
Since $P_q(q)=q^gP_q(1)$ and $q\not=\ell$, this implies 
\begin{equation}\label{ellP}
\ell \mid P_q(1).
\end{equation}
We first prove Lemma~\ref{Jac_p'}.
Let $\ell\nmid n$.
If $\ell\neq q_1$, taking $q=q_1$ in \eqref{ellP} gives
$\ell\mid P_{q_1}(1)$.
By condition~{\rm (i)}, this implies either $\ell\mid n$ or
$\ell=q_1$, a contradiction.
If $\ell=q_1$, taking $q=q_2$ gives
$q_1\mid P_{q_2}(1)$, contradicting condition~{\rm (ii)}.
Hence
\[
J(\mbb Q^{\mrm{ab}})[\ell]=0
\]
for every prime $\ell\nmid n$.
This proves Lemma~\ref{Jac_p'}.

Next we show Lemma~\ref{Jac_p}.
By \eqref{ellP}, 
we have  $X-1\mid P_{q_3}(X)$ in $\mbb{F}_{\ell}[X]$.
Since
$(X-1)^2\nmid P_{q_3}(X)$,
we have
$$
P_{q_3}'(1)\not\equiv0\pmod\ell.
$$
The functional equation
$$
q_3^gP_{q_3}(X)
=
X^{2g}P_{q_3}(q_3/X)
$$
also implies that
$X-q_3\mid P_{q_3}(X)$ modulo $\ell$.
Differentiating the functional equation and evaluating at $X=q_3$,
we obtain
$$
q_3^gP_{q_3}'(q_3)
\equiv
-q_3^{2g-1}P_{q_3}'(1)
\not\equiv0
\ \mrm{\mod} \ \ell.
$$
Hence 
$$
(X-q_3)^2\nmid P_{q_3}(X)
$$ in $\mbb{F}_{\ell}[X]$.
By \eqref{Caruso}, 
$\det(X-\mrm{Fr}_{q_3}\mid W_{\ell})
=\prod^d_{i=1}(X-q_3^{a_i})$
divides $P_{q_3}(X)$ mod $\ell$.
Since $a_i\in \{0,1\}$ and
$$
(X-1)^2\nmid P_{q_3}(X),
\qquad
(X-q_3)^2\nmid P_{q_3}(X)
$$
in $\mbb F_{\ell}[X]$, the number of indices $i$ with $a_i=0$ and the number
of indices $i$ with $a_i=1$ are both at most one.
Hence $d\le2$.
Since $\ell$ is an odd prime divisor of $n$, we have
$$
J(\mbb Q)[\ell]\cap
\Sigma(\overline{\mbb Q})[\ell]=0,
$$
because $G_{\mbb Q}$ acts trivially on the former
and via $\overline\chi_\ell$ on the latter. Thus the sum 
$J(\mbb Q)[\ell]+
\Sigma(\overline{\mbb Q})[\ell]$
is direct. Hence we have 
$J(\mbb Q)[\ell]\oplus
\Sigma(\overline{\mbb Q})[\ell]
\subset W_\ell$
and the left-hand side has dimension $2$ over $\mbb F_\ell$.
Therefore
$$
W_\ell
=
J(\mbb Q)[\ell]\oplus
\Sigma(\overline{\mbb Q})[\ell].
$$
This proves Lemma~\ref{Jac_p}.
\end{proof}

For any prime $q\neq N$, let $T_q$ denote the Hecke correspondence at $q$.
We use the same notation for the induced endomorphism of
$J=J_0(N)$ and the induced operator on $S_2(\Gamma_0(N))$.
For any prime $\ell\neq q,N$, we use the same notation
for the induced endomorphism of $T_\ell(J)$.

Let
$a_{q,1},\ldots,a_{q,g}$
be the eigenvalues of $T_q$ on
$S_2(\Gamma_0(N))$, counted with multiplicity, and put
$$
H_q(Y)
:=
\det\bigl(Y-T_q\mid S_2(\Gamma_0(N))\bigr)
=
\prod_{i=1}^g (Y-a_{q,i}).
$$
By the Eichler--Shimura relation
(see, for example, \cite[Chapter II, Section 6]{Maz77}), 
we have
$$
P_q(X)
=
\prod_{i=1}^g
(X^2-a_{q,i}X+q).
$$
Hence
\begin{equation}\label{PH}
P_q(X)
=
X^gH_q(X+q/X).
\end{equation}
Since the action of the Hecke algebra on
$S_2(\Gamma_0(N))$ is faithful,
the Cayley--Hamilton theorem gives
$$
H_q(T_q)=0
$$
as an equality in the Hecke algebra, and hence on $T_\ell(J)$.
Moreover,
$$
\det(X-T_q\mid T_\ell(J)\otimes_{\mbb Z_\ell}\mbb Q_\ell)
=H_q(X)^2.
$$
Indeed, each $a_{q,i}$ occurs twice
on $T_\ell(J)\otimes_{\mbb Z_\ell}\mbb Q_\ell$.

\begin{lemma}\label{Jac_p_infty}
Let $\ell$ be an odd prime divisor of $n$,
and put $m=v_{\ell}(n)$. 
Assume that there exist primes $q_3,q_4\notin \{\ell,N\}$
satisfying 
$(X-1)^2\nmid P_{q_3}(X)$
in $\mbb{F}_{\ell}[X]$
and 
$v_{\ell}(P_{q_4}(1))\le m$.
Then, 
$$
J(\mbb{Q}^{\mrm{ab}})[\ell^{\infty}]
=J(\mbb{Q}^{\mrm{ab}})[\ell^m]
=J(\mbb{Q})[\ell^m]\oplus \Sigma(\overline{\mbb{Q}})[\ell^m].
$$
\end{lemma}

\begin{proof}
First, we claim that 
\begin{equation}\label{Hq3}
(X-(q_3+1))\mid H_{q_3}(X) \quad  \mrm{and}\quad  
(X-(q_3+1))^2\nmid H_{q_3}(X)
\end{equation}
in $\mbb{F}_{\ell}[X]$.
Since $\ell\mid n$, we have
$J(\mbb Q)[\ell]\neq0$.
Since $J$ has good reduction at $q_3$
and $q_3\neq\ell$, 
the reduction map on $J(\mbb{Q})[\ell]$
at $q_3$ is injective.
Thus $\ell$ divides $\# J(\mbb{F}_{q_3})=P_{q_3}(1)=H_{q_3}(q_3+1)$. Thus
$(X-1)\mid P_{q_3}(X)$ and 
$(X-(q_3+1))\mid H_{q_3}(X)$ in $\mbb{F}_{\ell}[X]$.
Differentiating \eqref{PH} and evaluating at $X=1$, we obtain 
$$
P_{q_3}'(1)\equiv (1-q_3)H_{q_3}'(q_3+1) \ \mrm{mod}\ \ell,
$$
where we used
$H_{q_3}(q_3+1)\equiv0\pmod\ell$.
By assumption,
$(X-1)^2\nmid P_{q_3}(X)$
in $\mbb F_\ell[X]$.
Hence
$$
P_{q_3}'(1)\not\equiv 0 \ \mrm{mod}\ \ell.
$$
It follows that
$$
q_3\not\equiv1 \ \mrm{mod}\ \ell
\quad\text{and}\quad
H_{q_3}'(q_3+1)\not\equiv0 \ \mrm{mod}\ \ell.
$$
Therefore $q_3+1$ is a simple root of $H_{q_3}(X)$ modulo $\ell$,
which proves \eqref{Hq3}.

To simplify notation, we set
$T:=T_{\ell}(J)$ and
$W_r:=J(\mbb{Q}^{\mrm{ab}})[\ell^r]$.
By Lemma~\ref{Jac_p}, we have
$$
W_1
=
J(\mbb{Q})[\ell]
\oplus
\Sigma(\overline{\mbb{Q}})[\ell]
\simeq
(\mbb{Z}/\ell\mbb{Z})^2.
$$

We now construct a $G_{\mbb Q}$-equivariant idempotent
on $T$ which will isolate a rank-two direct summand
lifting $W_1$.
By \eqref{Hq3} and Hensel's lemma, we have a factorization
$$
H_{q_3}(X)=(X-a)h(X),
$$
where $a\in\mbb Z_\ell$ and
$h(X)\in\mbb Z_\ell[X]$ satisfy
$$
a\equiv q_3+1\pmod{\ell}
\quad\text{and}\quad
h(a)\in\mbb Z_\ell^\times.
$$
We define
$$
e:=\frac{h(T_{q_3})}{h(a)}
\in\mrm{End}_{\mbb Z_\ell}(T),
$$
and put
$$
U:=eT,
\qquad
U':=(1-e)T.
$$
For every $r\geq1$, we use the same notation $e$
for the endomorphism of $T/\ell^rT$ induced by $e$.

\begin{lemma}\label{idem}
The endomorphism $e$ has the following properties:
\begin{itemize}
\item[{\rm (1)}]
$e$ is $G_{\mbb Q}$-equivariant.
\item[{\rm (2)}]
$e^2=e$.
\item[{\rm (3)}]
$U=eT$ is a direct summand of $T$ of
$\mbb Z_\ell$-rank two.
\item[{\rm (4)}] For every $r\geq1$, the subgroup $W_r$ is stable under $e$.
Moreover, the induced endomorphism $e$ acts as the identity on $W_1$.
\item[{\rm (5)}] For every $r\geq1$, we have
$W_r=eW_r\subset U/\ell^rU$.
Moreover,
$$
W_r=J(\mbb{Q})[\ell^r]\oplus \Sigma(\overline{\mbb{Q}})[\ell^r]= U/\ell^rU
$$
for $1\leq r\leq m$.
\end{itemize}
\end{lemma}

\begin{proof}
Since $T_{q_3}$ commutes with the $G_{\mbb{Q}}$-action,
so does $e$. This shows (1).
Write $h(X)-h(a)=(X-a)g(X)$ for some $g(X)\in \mbb{Z}_{\ell}[X]$.
By the Cayley-Hamilton theorem, it follows that 
$e^2-e=h(a)^{-2}h(T_{q_3})(h(T_{q_3})-h(a))
=h(a)^{-2}g(T_{q_3})H_{q_3}(T_{q_3})=0$.
This shows (2).
By (2), we have a decomposition 
$$
T=U\oplus U',\ U:=eT,\ U':=(1-e)T.
$$
Since
$$
(T_{q_3}-a)e
=
h(a)^{-1}H_{q_3}(T_{q_3})
=0,
$$
the operator $T_{q_3}$ acts as multiplication by $a$ on $U$.
On the other hand,
$$
h(T_{q_3})(1-e)=h(a)e(1-e)=0,
$$
so $h(T_{q_3})=0$ on $U'$.
Since $h(a)\neq0$, the operator $T_{q_3}$ on
$U'\otimes_{\mbb Z_\ell}\mbb Q_\ell$
does not have $a$ as an eigenvalue.
Indeed, if $0\neq v\in U'\otimes\mbb Q_\ell$
satisfied $T_{q_3}v=av$, then
$0=h(T_{q_3})v=h(a)v$, a contradiction.
Therefore, using
$$
\det(X-T_{q_3}\mid T\otimes_{\mbb Z_\ell}\mbb Q_\ell)
=
H_{q_3}(X)^2
=
(X-a)^2h(X)^2,
$$
we obtain
$$
\mrm{rank}_{\mbb Z_\ell}U=2.
$$
This proves (3). 
Since $W_r=(T/\ell^rT)^{G_{\mbb{Q}^{\mrm{ab}}}}$
and $e$ is $G_{\mbb{Q}}$-equivariant,
$W_r$ is stable under $e$. 
By the Eichler-Shimura relation, we have 
$T_{q_3}=\mrm{Fr}_{q_3}+q_3\mrm{Fr}_{q_3}^{-1}$ on $T$.
Since $W_1=J(\mbb{Q})[\ell]\oplus \Sigma(\overline{\mbb{Q}})[\ell]
\simeq J(\mbb{Q})[\ell]\oplus \mu_{\ell}(\overline{\mbb{Q}})$,
$T_{q_3}$ acts as multiplication by $q_3+1$ on $W_1$.
Since
$a\equiv q_3+1\pmod{\ell}$,
we have
$h(q_3+1)\equiv h(a)\ \mrm{mod}\ \ell.$
As $h(a)\in\mbb Z_\ell^\times$, the endomorphism induced by $e$ on
$W_1\subset T/\ell T$ is given by
\[
h(a)^{-1}h(q_3+1)=1
\]
in $\mbb F_\ell$.
Therefore $e$ acts as the identity on $W_1$.
This proves (4). Finally, we show (5).
By (4), we have a decomposition 
$$
W_r=eW_r\oplus (1-e)W_r.
$$
Since $((1-e)W_r)[\ell]=(1-e)(W_r[\ell])=(1-e)W_1=0$, we obtain $(1-e)W_r=0$.
Thus $W_r=eW_r\subset U/\ell^rU$.
Suppose now that $r\le m$. Then 
$W_r\supset J(\mbb{Q})[\ell^r]\oplus \Sigma(\overline{\mbb{Q}})[\ell^r]\simeq (\mbb{Z}/\ell^r\mbb{Z})^2$, and the last term has the same order as
that of $U/\ell^rU$ by (3).
Therefore, we obtain 
$J(\mbb{Q})[\ell^r]\oplus \Sigma(\overline{\mbb{Q}})[\ell^r]=W_r=U/\ell^rU$
as desired.
\end{proof}

Now we return to the proof of Lemma \ref{Jac_p_infty}.
By Lemma \ref{idem} (5), for the proof of the lemma, 
it is enough to prove $J(\mbb{Q}^{\mrm{ab}})[\ell^{\infty}]=W_m$.
It suffices to prove 
$$
\ell^mW_{m+1}=0.
$$
Indeed, this implies $W_{m+1}=W_m$.
If there were a point of order $\ell^r$ in
$J(\mbb Q^{\mrm{ab}})$ for some 
$r\ge m+1$,
then an appropriate multiple of it would give a point
of order $\ell^{m+1}$, a contradiction.
Thus
$$
J(\mbb Q^{\mrm{ab}})[\ell^\infty]=W_m.
$$

Fix an embedding
$\overline{\mbb Q}\hookrightarrow\overline{\mbb Q}_{\ell}$.
We regard $I_{\ell}=I_{\mbb{Q}_{\ell}}$ and also 
regard $W_{m+1}$ as an
$I_{\mbb{Q}_{\ell}}$-stable submodule of
$$
J(\mbb Q_{\ell}^{\mrm{ab}})[\ell^{m+1}].
$$
Proposition~\ref{pro-nilpotent} shows that condition (C1) holds for
$J_{/\mbb Q_{\ell}}$ and $\mbb Q_{\ell}^{\mrm{ab}}/\mbb Q_{\ell}$.
Applying local decomposition  to the $I_{\ell}$-stable submodule
$W_{m+1}$ (cf.\ Theorem~\ref{LD:p-adic}), 
we obtain
$$
W_{m+1}=M_0\oplus M_1,\quad 
M_0:=W_{m+1}^{I_{\ell}},\quad  
M_1:=W_{m+1}\cap \hat{J}(\mbb{Q}^{\mrm{ab}}_{\ell})[{\ell}^{m+1}].
$$
By Corollary~\ref{hatG_ab}, the inertia group $I_{\ell}$ acts on $M_i$
via $\chi_{\ell}^i$ for $i=0,1$.
Choose $\tau_0\in I_{\ell}$ such that $\chi_{\ell}(\tau_0)$ is a topological
generator of $\chi_{\ell}(I_{\ell})=\mbb{Z}_{\ell}^{\times}$.
Then $M_i$ is characterized as the $\chi_{\ell}^i(\tau_0)$-eigenspace of
$\tau_0$ on $W_{m+1}$; namely,
$$
M_i=\{x\in W_{m+1}\mid \tau_0(x)=\chi_{\ell}^i(\tau_0)x\}.
$$
Since the action of $G_{\mbb Q}$ on $W_{m+1}$ factors through an
abelian quotient,  each $M_i$ is stable under the action of $G_{\mbb Q}$.
Since $W_{m+1}[{\ell}]=J(\mbb{Q})[\ell]\oplus \Sigma(\overline{\mbb{Q}})[\ell]
\simeq \mbb{F}_{\ell}\oplus \mbb{F}_{\ell}(1)$, we have 
$\dim_{\mbb{F}_{\ell}} M_i[{\ell}]=1$ for $i=0,1$.
Thus $M_i$ is cyclic.

\medskip 

Assume ${\ell}^mW_{m+1}\not=0$.
Then $M_i$ is cyclic of order ${\ell}^{m+1}$ for some $i$.
Let 
$$
\lambda_i\colon G_{\mbb{Q}}\to (\mbb{Z}/{\ell}^{m+1}\mbb{Z})^{\times}
$$
be the representation given by the $G_{\mbb{Q}}$-action on this $M_i$.

\begin{lemma}\label{lambda_ur}
$\lambda_i$ is unramified at $N$.
\end{lemma}
\begin{proof}
Since $J$ has semistable reduction at $N$ and ${\ell}\neq N$,
the action of $I_N$ on $T_{\ell}(J)$ is unipotent.
More precisely, Grothendieck's monodromy theorem gives that 
$\rho_U(\tau)-1$ is nilpotent 
for every $\tau\in I_N$, where
$$
\rho_U\colon G_{\mbb Q}\to GL_{\mbb Z_{\ell}}(U)
$$
denotes the representation induced by the $G_{\mbb Q}$-action on $U$.
Since $U/{\ell}^mU=W_m=J(\mbb{Q})[{\ell}^m]\oplus \Sigma(\overline{\mbb{Q}})[{\ell}^m]$ is unramified at $N$,  
we have 
$$
\rho_U(\tau)=1+{\ell}^mX_{\tau}
$$
for some nilpotent $X_{\tau}\in \mrm{End}_{\mbb{Z}_{\ell}}(U)$.
On the other hand,
since $M_i[{\ell}^m]\, (\subset W_m)$ is unramified at $N$,
we may write 
$$
\lambda_i(\tau)=1+{\ell}^mx_{\tau}
$$
for some $x_{\tau}\in \mbb{Z}/{\ell}^{m+1}\mbb{Z}$. 
By Lemma \ref{idem} (5), 
we can choose $v\in U$ whose image in $U/{\ell}^{m+1}U$ is a generator of $M_i$.
Then 
$$
\rho_U(\tau)v
\equiv
\lambda_i(\tau)v\ 
\mrm{mod}\ {{\ell}^{m+1}U},
$$
and hence
$$
X_\tau v\equiv x_\tau v \ \mrm{mod}\ {{\ell}U}.
$$
Since $v\notin {\ell}U$ and $X_{\tau}$ is nilpotent, 
we have $x_\tau\equiv0$ mod $\ell$.
Thus $\lambda_i(\tau)=1$ for every $\tau\in I_N$, proving that
$\lambda_i$ is unramified at $N$.
\end{proof}
We return to the proof of Lemma~\ref{Jac_p_infty}.
The character $\lambda_i\chi_{\ell}^{-i}$ is unramified at ${\ell}$
by the definition of $M_i$.
At every prime $q\neq {\ell},N$, it is unramified by the good reduction
of $J$, while it is unramified at $N$ by Lemma~\ref{lambda_ur}.
Hence $\lambda_i\chi_{\ell}^{-i}$ is unramified at every finite prime.
Since $\mbb Q$ has no nontrivial abelian extension unramified at every
finite prime, we obtain
$$
\lambda_i=\chi_{\ell}^i.
$$
Hence we have 
$\mrm{Fr}_qx=q^ix$
for any prime $q\notin\{\ell,N\}$ and 
any $x\in M_i$. 
Since $M_i$ is cyclic of order ${\ell}^{m+1}$ and 
$P_q(\mrm{Fr}_q)=0$ in $\mrm{End}_{\mbb{Z}_{\ell}}(T)$,
it follows that 
${\ell}^{m+1}\mid P_q(q^i)$. 
On the other hand, the functional equation
of the Frobenius polynomial gives 
$
P_q(q)=q^{\dim J}P_q(1).
$
Thus we have 
$$
{\ell}^{m+1}\mid P_q(1)
$$ 
for every prime $q\not={\ell},N$.
Taking $q=q_4$, we obtain
$$
v_\ell(P_{q_4}(1))\ge m+1,
$$
contradicting the assumption
$v_\ell(P_{q_4}(1))\le m$.
Therefore ${\ell}^mW_{m+1}=0$, as desired.
This completes the proof.
\end{proof}

\begin{theorem}\label{Jactor}
Assume that the assumptions of Lemma~\ref{Jac_p'} hold
and that, for every odd prime divisor $\ell$ of $n$,
the assumptions of Lemma~\ref{Jac_p_infty} hold.

\begin{itemize}
\item[{\rm (1)}]
If $n$ is odd, then
$$
J_0(N)(\mbb{Q}^{\mrm{ab}})_{\mrm{tor}}
=
J(\mbb{Q})_{\mrm{tor}}
\oplus
\Sigma(\overline{\mbb{Q}})_{\mrm{tor}}.
$$
In particular, this group is isomorphic to
$(\mbb{Z}/n\mbb{Z})^2$ as an abstract group.

\item[{\rm (2)}]
If $n$ is even, then
$$
J_0(N)(\mbb{Q}^{\mrm{ab}})_{2'}
=
J(\mbb{Q})_{2'}
\oplus
\Sigma(\overline{\mbb{Q}})_{2'}.
$$
In particular, this group is isomorphic to
$(\mbb{Z}/n_{\mrm{odd}}\mbb{Z})^2$
as an abstract group, where $n_{\mrm{odd}}$ denotes the odd part of $n$.
\end{itemize}
\end{theorem}

\begin{theorem}\label{Jacstructure}
Assume that $n$ is odd and that $J_0(N)$ is nontrivial.
Assume also that the assumptions of Lemma~\ref{Jac_p'} hold
and that, for every prime divisor $\ell$ of $n$,
the assumptions of Lemma~\ref{Jac_p_infty} hold.
Then
$$
J_0(N)(\mbb{Q}^{\mrm{ab}})
\simeq
\left(\bigoplus_{j=1}^{\infty}\mbb Z\right)
\oplus
(\mbb Z/n\mbb Z)^2.
$$
\end{theorem}

\begin{proof}
Since $J_0(N)(\mbb{Q}^{\mrm{ab}})_{\mrm{tor}}$ is finite,
it follows from \cite[Proposition 7]{Moo09} that
$$
J_0(N)(\mbb{Q}^{\mrm{ab}})
/
J_0(N)(\mbb{Q}^{\mrm{ab}})_{\mrm{tor}}
$$
is a free $\mbb Z$-module of at most countable rank.
On the other hand, $J_0(N)$ has a newform quotient $A_f$
of $\mrm{GL}_2$-type with trivial character.
By \cite[Theorem 1.3]{SaYa12},
$A_f(\mbb{Q}^{\mrm{ab}})$ has infinite rank.
Hence $J_0(N)(\mbb{Q}^{\mrm{ab}})$ also has infinite rank.
Therefore its torsion-free quotient is a free $\mbb Z$-module
of countably infinite rank.
Combining this with Theorem~\ref{Jactor}(1),
we obtain the desired result.
\end{proof}

A computation in MAGMA verifies the hypotheses of Theorem~\ref{Jactor} 
for every level $N$ listed in Table~\ref{tab:jacobian-torsion}.
We take $(q_1,q_2)$ as indicated in the table for Lemma~\ref{Jac_p'},
and $q_3=q_4=2$ in Lemma~\ref{Jac_p_infty} for every odd prime divisor
$\ell$ of $n$.

\begin{table}[htbp]
\centering
\small
\setlength{\tabcolsep}{5pt}
\begin{tabular}{@{}ccccc@{\qquad}ccccc@{}}
\hline
$N$ & $g$ & $n$ & $(q_1,q_2)$ & $J_0(N)(\mathbb Q^{\mathrm{ab}})_{\mathrm{tor}}$
& $N$ & $g$ & $n$ & $(q_1,q_2)$ & $J_0(N)(\mathbb Q^{\mathrm{ab}})_{\mathrm{tor}}$ \\
\hline
$11$ & $1$ & $5$  & $(2,3)$ & $C_5^2$
& $53$ & $4$ & $13$ & $(2,3)$ & $C_{13}^2$ \\
$19$ & $1$ & $3$  & $(2,5)$ & $C_3^2$
& $59$ & $5$ & $29$ & $(2,3)$ & $C_{29}^2$ \\
$23$ & $2$ & $11$ & $(2,3)$ & $C_{11}^2$
& $61$ & $4$ & $5$  & $(2,5)$ & $C_5^2$ \\
$29$ & $2$ & $7$  & $(2,3)$ & $C_7^2$
& $71$ & $6$ & $35$ & $(3,5)$ & $C_{35}^2$ \\
$37$ & $2$ & $3$  & $(7,2)$ & $C_3^2$
& $73^\dagger$ & $5$ & $6$ & $(19,3)$ & $C_3^2$ \\
$41^\dagger$ & $3$ & $10$ & $(3,2)$ & $C_5^2$
& $79$ & $6$ & $13$ & $(2,5)$ & $C_{13}^2$ \\
$47$ & $4$ & $23$ & $(2,3)$ & $C_{23}^2$
& $83$ & $7$ & $41$ & $(2,3)$ & $C_{41}^2$ \\
\hline
\end{tabular}
\caption{Torsion of modular Jacobians over $\mathbb Q^{\mathrm{ab}}$ 
with only the odd parts shown for daggered entries.}
\label{tab:jacobian-torsion}
\end{table}
Here $g=\dim J_0(N)$, $n=(N-1)/\gcd(N-1,12)$ and $C_m=\mathbb Z/m\mathbb Z$.
The following MAGMA code checks the hypotheses of
Lemmas~\ref{Jac_p'} 
and~\ref{Jac_p_infty} 
for primes $11\le N<100$ with odd $n>1$,
using $q_1,q_2\in\{2,3,5,7\}$ and $q_3=q_4=2$.

{\small
\begin{verbatim}
Z := Integers(); Q := [2,3,5,7]; L := [];
for N in PrimesInInterval(11,99) do
    n := (N-1) div GCD(N-1,12);
    if n eq 1 or IsEven(n) then continue; end if;
    J := JZero(N); P := FrobeniusPolynomial(J,2);
    a := Z!Evaluate(P,1);
    b := Z!Evaluate(Derivative(P),1);
    assert a mod n eq 0;
    if GCD(n,b) ne 1 or GCD(n,a div n) ne 1 then
        continue;
    end if;
    c := [Z!Evaluate(FrobeniusPolynomial(J,q),1) : q in Q];
    if exists{i : i in [1..#Q] |
        forall{l : l in PrimeDivisors(c[i]) |
            (n*Q[i]) mod l eq 0} and
        exists{j : j in [1..#Q] |
            j ne i and c[j] mod Q[i] ne 0}} then
        Append(~L,N);
    end if;
end for;
print L;
\end{verbatim}
}

The output is $[11,19,23,29,37,47,53,59,61,71,79,83]$.
\section{Torsion over relative ray class field towers}\label{sec:ray}
\subsection{Ray class towers}
The class field theory needed for the main theorem has two steps. 
We first describe the relative ray class group as a quotient of local units. 
We then use local-global compatibility of reciprocity to identify this quotient 
with the Galois group of the corresponding local extension. 

Throughout this section, fix a number field $F$, an odd rational prime $p$, 
a prime ideal $\mathfrak p \mid p$ of $F$, 
and a nonzero integral ideal $\mathfrak b$ 
relatively prime to $\mathfrak{p}$. 
Assume that $K:=F_{\mathfrak{p}}$ is unramified over $\mathbb{Q}_p$. 
Let $\mathfrak{m}_{\infty}$ be the product of all real places of $F$, 
with the empty product understood when $F$ has no real places. 
For $n \geq 0$, define $F_n/F$ to be the ray class field of modulus 
$\mathfrak{b} \mathfrak{p}^n \mathfrak{m}_{\infty}$ 
and put 
\[ H:=F_0,\quad  F_\infty:=\bigcup_{n \geq 0} F_n.\] 
The prime $\mathfrak{p}$ and the ideal $\mathfrak{b}$ remain fixed. 
We therefore suppress their dependence in the notation $F_n$. 
When it is useful to display the auxiliary modulus, the same field may be written 
$F_n(\mathfrak{b})$. 
The field $H/F$ is unramified at $\mathfrak{p}$. 
If $\mathfrak{b}=\mathcal{O}_F$, then $H=H^+$ is the narrow Hilbert class field of $F$. 

For a nonzero integral ideal $\mathfrak{c}$, 
let $\mathrm{Cl}_{\mathfrak{c}}^+$ be the ray class group of modulus $\mathfrak{c} \mathfrak{m}_{\infty}$. 
Thus $\mathrm{Cl}_{\mathfrak{c}}^+$ is the group of fractional ideals prime to $\mathfrak{c}$ modulo the principal 
ideals $(a)$ with $a \in 1+\mathfrak{c} \mathcal{O}_{F_{\mathfrak{q}}}$ for $\mathfrak{q} \mid \mathfrak{c}$ 
and $\sigma(a)>0$ for any real embedding $\sigma$. 
Let 
\[ E_{\mathfrak{c},1}^+:=\{\varepsilon \in \mathcal{O}_F^{\times} \mid \varepsilon \equiv 1 \pmod{\mathfrak{c}}, 
\, \sigma(\varepsilon)>0 \hspace{4pt} \text{for any real embedding} \hspace{4pt} \sigma \},\quad 
\mathcal{E}:=E_{\mathfrak{b},1}^+. \]
The positivity condition is vacuous when $F$ is totally imaginary. 
We view $\mathcal{E}$ as a subgroup of $\OO_K^{\times}$ through $F \hookrightarrow K$. 
Write $\Delta:F^{\times} \to \mathbb{A}_F^{\times}$ for the diagonal embedding 
and $\iota_{\pp}:K^{\times} \to \mathbb{A}_F^{\times}$ for the map whose $\pp$-component is the given element 
and whose other components are $1$. 
All Artin maps below use the same arithmetic Frobenius convention.
For any $p$-adic field $B$ occurring below, write
$U_B^{(0)}=\OO_B^\times$ and $U_B^{(n)}=1+\mathfrak m_B^n$ for
$n\ge1$, where $\mathfrak m_B$ is its maximal ideal.
We use $v_p(p)=1$ and $v_p(0)=\infty$. 

Choose a place $\mathfrak{P}$ of $F_\infty$ above $\mathfrak{p}$ and put 
\[
\mathfrak{P}_n:=\mathfrak{P}|_{F_n},\quad J:=H_{\mathfrak{P}_0},\quad M_n:=(F_n)_{\mathfrak{P}_n},\quad 
M_{\infty}:=\bigcup_{n \geq 0} M_n.
\]
We regard these fields, and the global fields embedded in them, as subfields of a fixed algebraic closure of $K$. 
Note that $M_0=J$, whereas $K=F_{\mathfrak{p}}$; these two local fields need not be equal. 
\begin{proposition}\label{prop:RayExactSeq}
For $n \geq 1$, put $\Gamma_n:=\mathrm{Gal}(F_n/H)$.
There is a natural exact sequence 
\[
 1\longrightarrow E^+_{\bb,1}/E^+_{\bb\pp^n,1}
 \xrightarrow{\ j_n\ }\OO_K^\times/\Up Kn
 \xrightarrow{\ \alpha_n\ }\Gamma_n\longrightarrow 1,
\]
where $j_n$ sends the class of a global unit to its local residue class, 
and $\alpha_n(u \Up Kn)=\Art_{F_n/F}(\iota_{\pp}(u))$. 
In particular, if $\widetilde{\alpha}_n:\OO_K^{\times} \to \Gamma_n$ is the map 
before taking the quotient by $\Up Kn$, then 
\begin{equation}\label{eq:ray-unit-kernel}
\ker \widetilde{\alpha}_n= \mathcal{E} \Up Kn, \quad 
\Gamma_n \simeq \OO_K^{\times}/(\mathcal{E} \Up Kn). 
\end{equation}
\end{proposition}
\begin{proof}
This is the specialization of the relative ray class sequence of Kopp--Lagarias~\cite[Theorem 6.5]{KL25} 
to the maximal order, comparing the ray class groups for the two ray moduli $\bb \pp^n \minfty$ and $\bb \minfty$. 
We give a direct proof for these moduli. 

For a nonzero integral ideal $\mathfrak{c}$, define the idele subgroup 
$U_{\mathfrak{c}}^+ \subset \mathbb{A}_F^{\times}$ by 
\[
(U_{\mathfrak{c}}^+)_{\mathfrak{q}}=\begin{cases} \OO_{F_\mathfrak{q}}^{\times}, & \mathfrak{q} \nmid \mathfrak{c} \\
1+\mathfrak{c} \OO_{F_{\mathfrak{q}}}, & \mathfrak{q} \mid \mathfrak{c} \end{cases}
\]
at finite places, with component $\mathbb{R}_{>0}$ at real places and 
$\mathbb{C}^{\times}$ at complex places. 
By the usual idelic description of ray class fields~\cite[Chapter VI]{Neu99}, 
global reciprocity gives $\mathrm{Gal}(F_n/F) \simeq \mathbb{A}_F^{\times}/(\Delta(F^{\times}) 
U_{\bb \pp^n}^+)$. 
The corresponding identification at $n=0$, together with restriction, yields
\[
\Gamma_n \simeq \frac{\Delta(F^{\times})U_{\mathfrak{b}}^+}{\Delta(F^{\times}) U_{\bb \pp^n}^+}.
\]
Since $\pp\nmid\bb$, the subgroups $U^+_{\bb}$ and
$U^+_{\bb\pp^n}$ differ only in their $\pp$-components. Hence
the embedding $\iota_{\pp}$ induces a surjective homomorphism
\[
 \OO_K^\times/\Up Kn
 \longrightarrow
 \frac{\Delta(F^\times)U^+_{\bb}}
      {\Delta(F^\times)U^+_{\bb\pp^n}}
 \simeq \Gamma_n.
\]
This is precisely the map $\alpha_n$.

It remains to compute its kernel. If $ u \in \OO_K^{\times}$ has trivial image, write 
\[ 
\iota_{\pp}(u)=\Delta(a)v,\qquad a \in F^{\times}, \quad v \in U_{\bb \pp^n}^+.
\]
Both $\iota_{\pp}(u)$ and $v$ belong to $U_{\bb}^+$, hence 
$\Delta(a) \in \Delta(F^{\times}) \cap U_{\bb}^+=\Delta(\mathcal{E})$. 
Indeed, the condition $\Delta(a)\in U^+_{\bb}$ means that
$a\in\OO_{F_{\mathfrak q}}^\times$ for every finite place 
$\mathfrak q\nmid\bb$, that
$a\in 1+\bb\OO_{F_{\mathfrak q}}$ for every
$\mathfrak q\mid\bb$, and that $a$ is positive at every real place.
Hence $a\in E^+_{\bb,1}$.
The $\pp$-component now gives $u \in \mathcal{E} \Up Kn$. 
Conversely, $\Up Kn$ maps trivially, and for $\varepsilon \in \mathcal{E}$ 
the idele $\iota_{\pp}(\varepsilon) \Delta(\varepsilon)^{-1}$ belongs to $U_{\bb \pp^n}^+$; 
its $\pp$-component is $1$, its components at primes dividing $\bb$ satisfy the 
required congruence conditions, 
and its real components are positive. Thus $\ker\widetilde{\alpha}_n=\mathcal{E} \Up Kn$. 
Since $\bb$ and $\pp^n$ are relatively prime, we have $\mathcal{E} \cap \Up Kn=E_{\bb \pp^n,1}^+$. 
\end{proof}
The compatibility with the usual ray class groups can also be recorded
as the commutative diagram
\[
\begin{CD}
 \Cl^+_{\bb\pp^n}(F) @>{\Art_n}>> \Gal(F_n/F)\\
 @V{\mathrm{pr}_n}VV @VV{\mathrm{res}}V\\
 \Cl^+_{\bb}(F) @>{\Art_0}>> \Gal(H/F).
\end{CD}
\]
Here $\mathrm{pr}_n$ forgets the congruence condition at $\pp^n$.
The group $\Gamma_n$ in Proposition~\ref{prop:RayExactSeq} is the kernel of the
right-hand vertical map. Thus the exact sequence concerns growth relative to
$H$, rather than the whole extension $F_n/F$.

\begin{corollary}\label{cor:ray-local}
The extension $J/K$ is unramified. 
For $n \geq 0$, the chosen prime $\mathfrak{P}_0$ of $H$ has a unique extension to $F_n$, 
this extension is totally ramified, and localization induces
\[
 \Gal(F_n/H)\simeq\Gal(M_n/J).
\]
In particular, $M_{\infty} \cap K^{\mathrm{ur}}=J$. 
Let $C$ be the closure of $\mathcal E$ in $\OO_K^\times$.
Then, for every $n\ge0$,
\[
 \Gal(F_\infty/F_n)
 \simeq \Gal(M_\infty/M_n)
 \simeq CU_K^{(n)}/C.
\]
\end{corollary}
\begin{proof}
The modulus of $H/F$ is prime to $\pp$, so $J/K$ is unramified.
For $n\ge1$, let
$\ell_n:\Gal(M_n/K)\hookrightarrow\Gal(F_n/F)$
be the embedding induced by localization.
Local--global compatibility of reciprocity
\cite[Chapter VI]{Neu99} gives
\begin{equation}\label{eq:local-global-diagram}
\begin{CD}
 K^\times @>{\iota_{\pp}}>> \mathbb A_F^\times\\
 @V{\Art_{M_n/K}}VV @VV{\Art_{F_n/F}}V\\
 \Gal(M_n/K) @>{\ell_n}>> \Gal(F_n/F).
\end{CD}
\end{equation}
Local reciprocity and Proposition~\ref{prop:RayExactSeq} give
\[
 I_{\mathfrak P_n}(F_n/F)
 =\ell_n\bigl(\Art_{M_n/K}(\OO_K^\times)\bigr)
 =\Gamma_n.
\]
Thus the relative inertia group
$I_{\mathfrak P_n}(F_n/H)=I_{\mathfrak P_n}(F_n/F)\cap\Gamma_n$
and the relative decomposition group both equal $\Gamma_n$.
Consequently $\mathfrak P_0$ has a unique extension to $F_n$, which
is totally ramified, and localization gives
$\Gal(F_n/H)\simeq\Gal(M_n/J)$.
The case $n=0$ is immediate. Since $J/K$ is unramified and each
$M_n/J$ is totally ramified, $M_\infty\cap K^{\mathrm{ur}}=J$.

By~\eqref{eq:ray-unit-kernel} and compactness of $\OO_K^\times$,
the finite-layer reciprocity maps induce a surjection
\[
 \OO_K^\times\longrightarrow\Gal(F_\infty/H)
\]
whose kernel is
\[
 \bigcap_{m\ge1}\mathcal E\Up Km=C.
\]
For every $n\ge0$, the preimage of $\Gal(F_\infty/F_n)$
under this map is $\mathcal E\Up Kn=C\Up Kn$,
since $\Up Kn$ is open.
Restricting this map and taking inverse limits of the
localization isomorphisms proves the last assertion.
\end{proof}
Uniqueness here is over the chosen prime $\mathfrak P_0$ of $H$.
There may be several primes of $H$ above $\pp$. Accordingly,
Corollary~\ref{cor:ray-local} compares the relative groups over $H$ and
$J$, and does not assert that the decomposition subgroup at $\mathfrak P_n$ 
equals $\Gal(F_n/F)$. 
This distinction will matter when the growth is
measured from $F$ rather than from $H$.

\subsection{Admissible uniformizers}

The next proposition compares the ray class layers with a single
Lubin--Tate tower. Put
\[
 \beta:=\inf_{\varepsilon\in\mathcal E}
 v_p\bigl(\Nr_{K/\mathbb Q_p}(\varepsilon)-1\bigr)
 \in\mathbb Z_{\ge0}\cup\{\infty\}.
\]
The norm is taken locally via the fixed embedding $F\hookrightarrow K$,
not globally over $F$; the quantity $\beta$ depends only on the tower.
In the remainder of this section, $0\le n\le\infty$ means
$n\in\mathbb Z_{\ge0}\cup\{\infty\}$.
We write $J_{\varpi,\infty}=J_{\varpi}$ and
$\mathfrak P_\infty=\mathfrak P$.
The invariant $b_{M_n}$ below is taken with $J$ as the base field.
\begin{proposition}\label{prop:D2}
Let $\varpi$ be any admissible uniformizer for $M_\infty/J$.
\begin{enumerate}
\item[(1)] For every $0\le n\le\infty$, we have
$M_n\subset J_{\varpi,n}$ and
\begin{equation}\label{eq:ray-intersection}
 F_\infty\cap J_{\varpi,n}=F_n.
\end{equation}

\item[(2)] For every $0\le n\le\infty$, we have
\[
 b_{M_n}=\min(\beta,n).
\]

\item[(3)] If $\pp=(\pi)$ has a totally positive generator and
$\pi^{[J:K]}\equiv1\pmod\bb$, then $\pi$ is admissible for
$M_\infty/J$, so one may take $\varpi=\pi$.
\end{enumerate}
\end{proposition}
\begin{proof}
For (1), first let $0\le n<\infty$.
Since $J/K$ is unramified, we have
$\Nr_{J/K}(\Up Jn)=\Up Kn$; see
\cite[Chapter V]{Neu99}.
By~\eqref{eq:local-global-diagram}, norm functoriality, and
Proposition~\ref{prop:RayExactSeq}, the group $\Up Jn$ acts
trivially on $M_n$ under reciprocity. For $n=0$, this follows
from $M_0=J$.
Together with admissibility, this gives
\[
 \varpi^{\mathbb Z}\Up Jn\subset\ker\Art_{M_n/J}.
\]
The left-hand side is the reciprocity kernel of $J_{\varpi,n}/J$,
so local class field theory gives $M_n\subset J_{\varpi,n}$.

Thus $F_n\subset F_\infty\cap J_{\varpi,n}$.
For the reverse inclusion, take $x\in F_\infty\cap J_{\varpi,n}$
and put $E=H(x)\subset F_m$ for some $m$.
Then $E/F$ is abelian, and the fixed embedding induces a place
$\mathfrak Q$ such that
$J\subset E_{\mathfrak Q}\subset J_{\varpi,n}$.
Local class field theory and norm transitivity give
\[
 \Up Kn=\Nr_{J/K}(\Up Jn)
 \subset\Nr_{E_{\mathfrak Q}/K}(E_{\mathfrak Q}^{\times}).
\]
Thus the conductor exponent of $E/F$ at $\pp$ is at most $n$.
Together with $E\subset F_m$, this shows that the 
conductor of $E/F$ divides 
$\mathfrak b\pp^n\mathfrak m_\infty$.
The conductor theorem \cite[Chapter VI]{Neu99} gives 
$E\subset F_n$, proving~\eqref{eq:ray-intersection}.
The assertions for $n=\infty$ follow by taking unions.

For (2), put
\[
 \mathcal U_n:=\ker\bigl(\OO_J^\times
 \xrightarrow{\Art_J}\Gal(M_n/J)\bigr)
 \qquad(0\le n\le\infty).
\]
By~\eqref{eq:local-global-diagram}, norm functoriality, and
Corollary~\ref{cor:ray-local}, we have
\[
 \mathcal U_\infty
 =\{u\in\OO_J^\times:\Nr_{J/K}(u)\in C\}.
\]
Since the norm on units is surjective,
$\Nr_{J/K}(\mathcal U_\infty)=C$.
For each $j\ge0$, the congruence
$\Nr_{K/\mathbb Q_p}(c)\equiv1\pmod{p^j}$ holds for every $c\in C$ if and only if it holds for every
$c\in\mathcal E$, by density and continuity.
Norm transitivity therefore gives $b_{M_\infty}=\beta$.

For $1\le n<\infty$, \eqref{eq:ray-unit-kernel} similarly gives
\[
 \mathcal U_n
 =\{u\in\OO_J^\times:\Nr_{J/K}(u)\in\mathcal E\Up Kn\},
 \qquad
 \Nr_{J/K}(\mathcal U_n)=\mathcal E\Up Kn.
\]
As $K/\mathbb Q_p$ is unramified, the absolute norm image is
\[
 \Nr_{J/\mathbb Q_p}(\mathcal U_n)
 =\Nr_{K/\mathbb Q_p}(\mathcal E)(1+p^n\mathbb Z_p).
\]
The ideal generated by $ab-1$, with
$a\in\Nr_{K/\mathbb Q_p}(\mathcal E)$ and
$b\in1+p^n\mathbb Z_p$, is $p^{\min(\beta,n)}\mathbb Z_p$:
setting $a=1$ or $b=1$ gives one inclusion, and
$ab-1=a(b-1)+(a-1)$ gives the other.
Thus $b_{M_n}=\min(\beta,n)$.
For $n=0$, the absolute norm image of
$\mathcal U_0=\OO_J^\times$ is $\mathbb Z_p^\times$,
so $b_{M_0}=0$ because $p>2$.

For (3), put $d=[J:K]$.
Then $\pi$ is a uniformizer of $J$ and
$\Nr_{J/K}(\pi)=\pi^d$.
The positivity and congruence assumptions imply
\[
 \iota_{\pp}(\pi^d)\Delta(\pi^d)^{-1}
 \in U^+_{\bb\pp^n}\qquad(0\le n<\infty).
\]
Global reciprocity and~\eqref{eq:local-global-diagram} give
$\Art_{M_n/K}(\pi^d)=1$.
By norm functoriality, $\Art_{M_n/J}(\pi)=1$ for every finite $n$,
so $\pi$ is admissible for $M_\infty/J$.
\end{proof}

\subsection{Decomposition and descent to finite layers}\label{subsec:descent}
Let $A/F$ be an abelian variety. 
Write $A(E)_{p'}$ for the prime-to-$p$ torsion subgroup of $A(E)$,
put $\Gamma_n=\Gal(F_n/H)$, and define
\[
 D_n:=D_{A,H}(F_n)
 :=\sum_{\sigma\in\Gamma_n}(\sigma-1)A(F_n)[p^\infty]
 \qquad(0\le n\le\infty).
\]
\begin{corollary}\label{cor:ray-total-torsion}
Let $A/F$ be an abelian variety which acquires semistable reduction 
over a finite extension $K'/K$ with $e(K'/K)<p-1$.
Then
\begin{equation}\label{eq:global-p-decomposition}
 A(F_n)_{\tor}
 =D_n\oplus A(H)[p^\infty]
   \oplus A(F_{\min\{n,1\}})_{p'}
 \qquad(0\le n\le\infty).
\end{equation}
\end{corollary}
\begin{proof}
By Corollary~\ref{cor:ray-local}, $J/K$ is unramified,
$M_n/J$ is totally ramified, and 
$\Gamma_n\simeq\Gal(M_n/J)$.
Moreover, $e(K'J/J)=e(K'/K)<p-1$ and $e_J=1$,
so Theorem~\ref{thm:intro-local} applies over $J$.
For each integer $n\ge1$, put $N_n=A(F_n)[p^\infty]$ 
and let $\Delta_n$ be the subgroup of $\Gamma_n$ 
consisting of all elements of order prime to $p$.
The projection
\[
 e_n:=|\Delta_n|^{-1}\sum_{\delta\in\Delta_n}\delta
 \colon N_n\longrightarrow N_n
\]
has kernel $D_n$ and image 
$N_n^{\Delta_n}=N_n^{\Gamma_n}=A(H)[p^\infty]$.
By Proposition~\ref{prop:RayExactSeq}, $F_1/H$ has degree prime 
to $p$ and $\Gal(F_n/F_1)$ is a $p$-group. 
Thus $M_n\cap J^{\mathrm{tame}}=M_1$, and the prime-to-$p$ part of 
the local decomposition gives $A(F_n)_{p'}=A(F_1)_{p'}$.
This proves the assertion for each integer $n\geq1$.
The case $n=0$ is immediate, and passage to the union proves the case $n=\infty$.
\end{proof}
Suppose now that $A$ has semistable reduction at $\pp$.
Put $A_J:=A\times_F J$ and $g=\dim A$. 
Write
\[
 0\longrightarrow T\longrightarrow G_J\longrightarrow B
 \longrightarrow0,\qquad
 G_J^{\mathrm{an}}\xrightarrow{u}A_J^{\mathrm{an}}
\]
for the Raynaud uniformization over $J$, with notation as in 
Section~\ref{LD_st}. Let $\mathcal G_J/\OO_J$ be the semi-abelian Raynaud model of $G_J$.
Propositions~\ref{LD:semistable} and~\ref{pro-nilpotent} give
\begin{equation}\label{eq:ray-rigid-intersection}
 D_n=A(F_n)[p^\infty]\cap
 u\bigl(\widehat G_J(M_n)[p^\infty]\bigr).
\end{equation}
Indeed, the local direct-sum decomposition restricts to the global 
Galois-stable subgroup, and its projection onto the ramified summand 
is $1-e_n$ at each finite layer.
The map $u$ is injective on torsion, so localization identifies $D_n$ 
with a subgroup of $\widehat G_J(M_n)[p^\infty]$.

If $A$ has good reduction at $\pp$, then $G_J=A_J$ and 
\eqref{eq:ray-rigid-intersection} becomes
\[
 D_n=\widehat A_{\mathfrak P_n}(F_n)[p^\infty],\qquad
 \widehat A_{\mathfrak P_n}(F_n)
 :=\ker\bigl(A(F_n)\longrightarrow\overline A_J(k_J)\bigr).
\]
For semistable reduction we retain $D_n$; the relevant local group is the image under $u$ of
$\widehat G_J(M_n)[p^\infty]$.

We use Fontaine's bound for finite flat group schemes to descend the $p$-primary torsion. 
For a finite Galois extension $E/J$, 
let $G^v$ denote the upper ramification groups in Serre's convention,
and put
\[
 G^{(v)}=G^{v-1}\quad(v\ge0),\qquad
 u(E/J)=\sup\{v\ge0:G^{(v)}\ne1\},
\]
with $u(E/J)=0$ for the trivial extension.
For an abelian extension, the Hasse--Arf theorem \cite[Chapter IV]{Ser68} and local reciprocity
\cite[Chapter V]{Neu99} give
\[
 c(E/J):=\min\{j\ge0:U_J^{(j)}\subset\Nr_{E/J}(E^\times)\}
 =u(E/J)\in\mathbb Z_{\ge0}.
\]
Fontaine's theorem~\cite[\S2, Th\'eor\`eme 1]{Fon85} states that 
the splitting field of a finite flat commutative group scheme 
$\mathscr H/\OO_J$ killed by $p^r$ satisfies
\begin{equation}\label{eq:Fontaine-bound}
 u(J(\mathscr H)/J)\le e_J\left(r+\frac1{p-1}\right).
\end{equation}
Since $e_J=1$ and $p>2$, every abelian subextension of this 
splitting field has conductor exponent at most $r$.
\begin{proposition}\label{prop:ray-Fontaine}
If $A/F$ has semistable reduction at $\pp$, then
\[ 
A(F_\infty)[p^r]=A(F_r)[p^r] \qquad (r\ge1).
\]
\end{proposition}
\begin{proof}
Fix an admissible uniformizer $\varpi$ for $M_\infty/J$. 
Let $P \in D_\infty[p^r]$. 
By~\eqref{eq:ray-rigid-intersection} and Lemma~\ref{hatG1}, 
$P=u(Q)$ for a unique $Q \in \widehat G_J(M_\infty)[p^r]$. 
The injectivity and Galois equivariance of $u$ on $p^r$-torsion give $J(P)=J(Q)$.
The point $Q$ belongs to the generic fiber of the finite flat group scheme $\mathcal G_J[p^r]/\OO_J$.
Thus $J(P)/J$, which is abelian, has conductor exponent at most $r$ by~\eqref{eq:Fontaine-bound}.
Together with the admissibility of $\varpi$, 
this gives $\varpi^{\mathbb Z}U_J^{(r)}
 \subset \ker\operatorname{Art}_{J(P)/J}$,
so $J(P)\subset J_{\varpi,r}$.
Proposition~\ref{prop:D2}(1) now implies $P\in A(F_\infty\cap J_{\varpi,r})=A(F_r)$.
Taking $p^r$-torsion in~\eqref{eq:global-p-decomposition}, we obtain
\[
 A(F_\infty)[p^r]
 =D_\infty[p^r]\oplus A(H)[p^r]\subset A(F_r)[p^r].
\]
The reverse inclusion is immediate.
\end{proof}

\subsection{Smith invariants and torsion growth}
Assume that $A/F$ has semistable reduction at $\pp$, 
and retain the Raynaud uniformization over $J$.
Put
\[
 t=\rank T,\quad f_B=\text{the $p$-rank of }\overline B,\quad
 h=t+f_B,\quad
 \Lambda=T_p(\mathcal M_{G_J}),\qquad
 \mathcal M_{G_J,n}=(\mathcal G_J[p^n]^D)^{\aet}.
\]
The module $\Lambda$ is free of rank $h$.
Choose an admissible uniformizer $\varpi$ for $M_\infty/J$, 
and set
\[
 s=\Nr_{J/\mathbb Q_p}(\varpi),\qquad
 q_J=\#k_J=p^{f_J},\qquad c_\varpi=q_J/s.
\]
Following Section~\ref{FG}, 
write $\Fr_{G_J}=\rho_{G_J}(\Art_J(\varpi))$ for arithmetic Frobenius on $\Lambda$, 
and use the same notation for its matrix in a chosen $\mathbb Z_p$-basis.
It is independent of the chosen uniformizer because $\rho_{G_J}$ is unramified.
Let
\[
 0\le a_1\le\cdots\le a_h\le\infty
\]
be the Smith exponents of $c_\varpi I_h-\Fr_{G_J}$, 
where a zero diagonal entry has exponent $\infty$. 
We use $C_{p^\infty}=\mathbb Q_p/\mathbb Z_p$ as in Section~\ref{FG}. 
An empty direct sum is zero, and the determinant of the empty matrix is $1$.

The descriptions of $\widehat G_J(M_n)[p^\infty]$ 
below require condition {\rm (C2)} 
for $A_J$ and $M_\infty/J$.
By the proof of Proposition \ref{C2} (c), condition (C2) holds if $B$ has ordinary reduction.
It also holds if the tame ramification index of $M_\infty/J$ divides $p-1$, 
or if $\gcd([J:\mathbb Q_p],(2\dim B)!)=1$; 
the latter two criteria follow by applying Proposition~\ref{C2} to $B$.
In particular, if $K=\mathbb Q_p$, Corollary~\ref{cor:ray-local} 
identifies $\operatorname{Gal}(M_\infty/J)$ with a quotient of $\mathbb Z_p^\times$, 
so its tame ramification index divides $p-1$.

\begin{theorem}\label{un4:thm:G}
For {\rm (2)}--{\rm (4)},
assume in addition that {\rm (C2)} holds for $A_J$ and
$M_\infty/J$.
\begin{enumerate}
\item[(1)] If $A$ has good reduction at $\pp$, then,
for every $0\le n\le\infty$,
\begin{equation}\label{un4:eq:G-decomposition}
A(F_n)_{\tor}
=\widehat A_{\mathfrak P_n}(F_n)[p^\infty]
\oplus A(H)_{\tor}.
\end{equation}

\item[(2)] For every $0\le n\le\infty$,
localization and uniformization give an injection
\begin{equation}\label{un4:eq:G-infty}
 D_n\hookrightarrow\widehat G_J(M_n)[p^\infty]
 \simeq\bigoplus_{i=1}^h C_{p^{\min(a_i,\beta,n)}},
\end{equation}
where the isomorphism is non-canonical.

\item[(3)] Put
\[
 R=\sum_{i=1}^h\min(a_i,\beta),\qquad
 c=\max_{1\le i\le h}\min(a_i,\beta),
\]
with $R=c=0$ when $h=0$.
If $c<\infty$, then $A(F_\infty)_{\tor}$ is finite and
\begin{equation}\label{un4:eq:G-bound}
 [A(F_\infty)[p^\infty]:A(H)[p^\infty]]
 =\#D_\infty\le p^R.
\end{equation}
Moreover,
\begin{equation}\label{un4:eq:G-stable}
 A(F_\infty)[p^\infty]=A(F_c)[p^\infty],\qquad
 A(F_\infty)_{\tor}=A(F_{\max(c,1)})_{\tor}.
\end{equation}
If $A$ has good reduction, the full torsion stabilizes at $F_c$,
and the index in~\eqref{un4:eq:G-bound} equals 
$[A(F_\infty)_{\tor}:A(H)_{\tor}]$.
In particular, $\beta<\infty$ implies $c\le\beta$ and 
$R\le h\beta$, without any Frobenius nonvanishing assumption.

\item[(4)] Suppose $\det(c_\varpi I_h-\Fr_{G_J})\ne0$.
Let
$X^*(T):=\operatorname{Hom}_{\overline J}
(T_{\overline J},\mathbb G_m)$ 
be the character group of $T$, and let
$\mathrm{Fr}_{X^*(T)}$ denote the action of
arithmetic Frobenius on $X^*(T)\otimes\mathbb Z_p$.
Let $P_{B,J}$ be the Frobenius polynomial of $\overline B$, 
with $P_{B,J}=1$ if $B=0$.
Then
\begin{align}\label{eq:ray-semistable-determinant}
 m:=\sum_{i=1}^h a_i
 &=v_p\det(c_\varpi I_h-\Fr_{G_J})\notag\\
 &=v_p(P_{B,J}(s))-f_J\dim B
   +v_p\det(c_\varpi I_t-\Fr_{X^*(T)}).
\end{align}
Consequently $R\le m$ and $c\le m$.
If also $\beta=\infty$, then $R=m$.
In the good reduction case, the nonvanishing hypothesis 
is equivalent to $P_{A,J}(s)\ne0$, and
\[
 m=v_p(P_{A,J}(s))-g f_J,\qquad
 A(F_\infty)_{\tor}=A(F_c)_{\tor}=A(F_m)_{\tor}.
\]
\end{enumerate}
\end{theorem}
\begin{proof}
For (1), assume that $A$ has good reduction at $\pp$.
Equation~\eqref{eq:ray-rigid-intersection} gives 
$D_n=\widehat A_{\mathfrak P_n}(F_n)[p^\infty]$.
Since $M_n/J$ is totally ramified, 
the N\'eron--Ogg--Shafarevich criterion and the Galois equivariance of localization give $A(F_n)_{p'}=A(H)_{p'}$.
The assertion follows from Corollary~\ref{cor:ray-total-torsion}.

For (2), taking
$\operatorname{Gal}(M_\infty/M_n)$-invariants
in the identifications of \textup{(C2)} shows that
\textup{(C2)} also holds for $A_J$ and $M_n/J$.
The uniformizer $\varpi$ is admissible for $M_n/J$,
and Proposition~\ref{prop:D2}(2) gives $b_{M_n}=\min(\beta,n)$.
The asserted isomorphism therefore follows from 
Theorem~\ref{FG:Main}(1), while the injection follows 
from~\eqref{eq:ray-rigid-intersection}.

For (3), assume $c<\infty$. 
The group on the right of 
\eqref{un4:eq:G-infty} at $n=\infty$ has order $p^R$ and exponent $p^c$.
This proves the bound in (3).
For $c\ge1$, Proposition~\ref{prop:ray-Fontaine} places 
$D_\infty$ in $A(F_c)[p^c]$. 
If $c=0$, then $D_\infty=0$.
The decomposition~\eqref{eq:global-p-decomposition} then proves 
both stabilization assertions and finiteness.
For good reduction, the prime-to-$p$ summand is already defined 
over $H$, which gives the sharper assertion.

For (4), the exact sequence in Theorem~\ref{hatG}, 
Cartier duality, and passage to Tate modules give an exact 
sequence of unramified $\mathbb Z_p$-representations
\[
 0\longrightarrow T_p(\overline B^\vee)
 \longrightarrow\Lambda
 \longrightarrow X^*(T)\otimes\mathbb Z_p\longrightarrow0.
\]
The determinant therefore factors into its abelian and toric parts. 
For the abelian part, the unit-root computation in the 
proof of Theorem~\ref{FG:Main}(2), applied to $B$, gives
\[
 v_p\det(c_\varpi-\Fr_B)
 =v_p(P_{B,J}(s))-f_J\dim B.
\]
This proves~\eqref{eq:ray-semistable-determinant}.
When $A$ has good reduction, 
$T=0$ and $B=A_J$; the remaining claims follow from (3) and $c\le m$.
\end{proof}
\begin{corollary}\label{cor:ray-finite-layer}
If $A/F$ has semistable reduction at $\pp$ and {\rm (C2)} holds
for $A_J$ and $M_\infty/J$, then
\[
 D_n=D_\infty[p^n]\qquad(0\le n\le\infty)
\]
and
\begin{equation}\label{un4:eq:G-level}
 A(F_n)_{\tor}
 =D_\infty[p^n]\oplus A(H)[p^\infty]
  \oplus A(F_{\min(n,1)})_{p'}.
\end{equation}
For good reduction this is 
$A(F_n)_{\tor}=\widehat A_{\mathfrak P}(F_\infty)[p^n] \oplus A(H)_{\tor}$.
\end{corollary}
\begin{proof}
The cases $n=0$ and $n=\infty$ are immediate.
Let $n\ge1$ be an integer.
Theorem~\ref{un4:thm:G}(2) gives $p^nD_n=0$, so
$D_n\subset D_\infty[p^n]$.

Conversely, Proposition~\ref{prop:ray-Fontaine} gives
$D_\infty[p^n]\subset A(F_n)[p^n]$.
For $P\in D_\infty[p^n]$, its unique lift in
$\widehat G_J(M_\infty)[p^n]$ is fixed by
$\operatorname{Gal}(M_\infty/M_n)$ and hence belongs to
$\widehat G_J(M_n)[p^n]$.
Equation~\eqref{eq:ray-rigid-intersection} therefore gives $P\in D_n$.
The decomposition follows from Corollary~\ref{cor:ray-total-torsion}.
\end{proof}
\begin{proof}[{Proof of Theorem \ref{thm:intro-ray} }]
Take $\bb=\mathcal O_F$, so that $H=H^+$. 
Since $\pp=(\pi)$ with $\pi$ totally positive, the ideal $\pp$ has trivial class 
in the narrow class group and hence splits completely in $H^+/F$. 
Thus $J=K$. 
Proposition~\ref{prop:D2}(3) shows that $\pi$ is admissible for $M_\infty/K$, 
and condition (C2) holds by Proposition \ref{C2}(c).
The bounds and stabilization statement now follow from 
Theorem~\ref{un4:thm:G}(3)--(4), since $f_J=[K:\Q_p]$.
If $s=q=\#k_K$, the functional equation gives 
$P_{A,K}(q)=q^g P_{A,K}(1)=q^g \#\overline{A}(k_K)$,
proving the final assertion.
\end{proof}

\bibliographystyle{alpha}
\bibliography{references}

\end{document}